\documentclass[11pt,twoside,a4paper]{article}

\usepackage[american]{babel}
\usepackage{bm}
\usepackage{amsfonts}
\usepackage{graphicx}
\usepackage{amsmath}
\usepackage{fancyhdr}
\usepackage{a4wide}
\usepackage{comment}
\usepackage{setspace}

\usepackage[utf8]{inputenc}

\usepackage{amsthm}

\usepackage[hidelinks]{hyperref}
\usepackage{multirow}
\usepackage{tikz}
\usepackage[ruled,vlined]{algorithm2e}
\usepackage{enumitem}
\usepackage{float}

\newtheorem{theorem}{Theorem}[section]

\newtheorem{lemma}{Lemma}[section]

\newtheorem{corollary}{Corollary}[section]

\graphicspath{ {figure/} }

\newcommand{\nn}{\nonumber}

\newcommand{\disp}{\displaystyle}

\usepackage{color}

\newcommand{\zd}{\delta}

\newcommand{\zs}{\sigma}

\newcommand{\ze}{\varepsilon}

\newcommand{\zg}{\gamma}

\newcommand{\zm}{\mu}

\newcommand{\zr}{\rho}

\newcommand{\zy}{\psi}
\newcommand{\zf}{\varphi}

\newcommand{\zh}{\eta}

\newcommand{\zO}{\Omega}

\newcommand{\dif}{\; \textrm d }

\newcommand{\dsp}[2]{\frac{  \partial^2 #1 }{  \partial #2^2   } }
\newcommand{\dst}[2]{\frac{  \dif^2 #1 }{  \dif #2^2   } }
\newcommand{\dpp}[2]{\frac{  \partial #1 }{  \partial #2   } }
\newcommand{\dpt}[2]{\frac{  \dif #1 }{  \dif #2   } }

\theoremstyle{definition}

\begin{document}
	
	\begin{center}
		\Large{\textbf{Deterministic classical limit of the optimal control problem of quantum particles with spin.}}
		\\[1cm]
		\small{O. Morandi}	\\
		\vskip0.5cm
		\textit{\textsf University of Florence  \\
			Firenze, Italy}
		\vskip0.5cm
		\textit{omar.morandi@unifi.it}
	\end{center}

	\begin{center}
		\begin{minipage}[h]{0.8\textwidth}
			\section*{\textsf{Abstract}}
			\small
			\textsf{}
			We study the optimal control problem applied to a gas of particles with spin confined in a material with Rashba spin-orbit coupling effect, in the presence of an external magnetic field. The evolution of the particle gas is described in the Wigner formalism. We investigate the classical limit of the optimal control problem, and we prove the convergence of the solution of the quantum problem toward the solution of a simplified optimal control model, based on an ODE description of the particle gas in terms of a single spin vector traveling along a classical trajectory in the phase-space.  
		\end{minipage}
	\end{center}
	\vspace{1cm}
	\normalsize

	\section{Introduction}
	
	Steering a quantum system towards a target state is of primary relevance in quantum information science \cite{Nielsen_10}. Today, control protocols applied to quantum systems are emerging as fundamental aspects for designing fault-tolerant experimental setups and engineering robust quantum processes. Optimal control algorithms find applications in various modern research areas such as, for example, the development of computers or communication devices applied to quantum information technology and the application of quantum entanglement in sensing and metrology \cite{Koch_22,Saffman_19,Henriet_20,Bluvstein_23}. 
	Optimal transport has also found several applications to ultra-cold atoms in optical traps \cite{Couvert_08,Murphy_09,Chen_11,Cicali_24,Torrontegui_11,Negretti_13,GueryOdelin_19}, including Bose-Einstein condensates \cite{Rosi_13,Jaeger_14}. Optimized control protocols are now available to enhance transport fidelity and suppress unwanted excitations \cite{Buecker_13}. Optimal control techniques have been proposed to achieve fast, high-fidelity transport of quantum systems, ensuring minimal energy cost and robustness against experimental imperfections \cite{Hohenester_07,Pagano_24,Muller_22}.
	Controlling confined or nearly free gases of charges particles subject to strong spin-orbit interaction is a relevant issue for spintronics. In spintronic devices, information is encoded in electron spin degrees of freedom, and spin-orbit interaction may be exploited to manipulate the spin orientation by electrical signals.
	Optimal Control (OC) procedures can be applied to design and engineer spintronic devices aimed at performing specific tasks such as, for example, generating continuous or modulated current of spin or modifying the spin orientation of localized or moving charges.  
	In the context of the mathematical models for optimal control of electron gases applied to spintronics, we investigate the well-posedness of an optimization problem where the evolution of a quantum gas with spin is controlled in the presence of the Zeeman and Rashba interaction. Our results give particular emphasis to the classical limit of such an optimal control problem and to the coherence of the latter with the full quantum dynamics.
	Optimal control aims to design a set of time-dependent parameters modulating the total electric potential with two main objectives: maximize the probability of finding the particle at the final time in a quantum state that may be considered as close as possible to a given classical phase-space coordinate (we will make this statement more precise in the following), and mean spin polarization aligned to a chosen direction. At the same time, the control should require the minimum possible source of external energy. 
	
	We study the classical limit of the optimal control problem. At this aim, we introduce a suitable $\hbar$ scaling of the initial datum of the system, in which the quantum dynamics tends to the deterministic classical transport of a single trajectory in phase-space and the spin vector rotates around to the total magnetic field given by the sum of the external and the Rashba effective fields. 
	The initial state of the quantum gas is described by a coherent state, and the associated Wigner distribution admits limit as a measure in phase-space, which can be characterized by a Dirac's delta. 

	Ensuring compatibility between classical dynamics and quantum evolution of a particle in the limit where the particle action is large compared to the Planck constant $\hbar$, has already been of primary concern in the early development of quantum mechanics. Since then, the classical limit of a quantum system has been widely investigated from both a physical and a mathematical point of view. 
	With the aim of providing solid mathematical ground to the description of the classical limit of a particle gas, around 80' it was recognized that the Wigner phase-space representation of quantum dynamics may constitute an ideal framework to investigate the classical limit. 
	Passing from the micro- to the macro- scale, the Wigner function ``behaves as expected" as a particle distribution. This can be seen by considering, as an example, the $\hbar\rightarrow0$ limit of a minimum uncertainly wave packet. Such a coherent state is represented in the Wigner description by a Gaussian function in the phase-space, whose position and momentum variances saturate the Heisenberg relations. Letting formally the Planck constant $\hbar$ go to zero, the Gaussian profile shrinks to a point, while its maximum diverges as $\hbar^{-\frac{3}{2}}$. This behavior clearly indicates that the convergence of the Wigner distribution is weak and should be found in the distribution space. The precise characterization of the classical limit of the Wigner pseudo-distribution function in terms of a phase-space measure was provided by Lions and Paul in \cite{Lions_93} and in its seminal contributions by Gérard \cite{Gerard_91}. The mathematical investigation of the classical limit of the Wigner distribution for linear and nonlinear systems was investigated by several authors \cite{Markowich_89,Markowich_Ringhofer,Panati_02,ArnoldRinghofer1996,ArnoldRinghofer1996,Pulvirenti_06,Gat_14,Denittis_03,Wu_12}. The concept of the Wigner measure has been fruitful in inspiring various generalizations such as, for example, the definition of the Bohm measure \cite{Sparber_03,Markowich_10,Markowich_12,Morandi_22_PLA,Figalli_14} associated with the Bohm formulation of quantum mechanics, and localisation principles for one-scale measure spaces \cite{Antonica_17}.
	
	Wigner formalism has been applied only recently to optimal control \cite{Morandi_PRA_25,Morandi_SIAM_25}, so the classical limit in this framework is unexplored and little is known regarding the classical limit of the problem of optimal control of quantum systems. The study of a classical limit in optimal control has dual interests. From a theoretical point of view the interest relies on the study of asymptotics of a specific class of nonlinear pseudo-differential equations associated to the quantum dynamics of a gas with spin. Concerning the applications, the classical limit is useful for deriving approximated models, which simplify the complex quantum dynamics and speed up the numerical computation associated with the optimization process.  
	
	\section{Quantum dynamics}\label{Sec_model}
	
	We consider the dynamics of a quantum particle gas with spin in the presence of the Rashba field, an external electric potential energy $U$ and a magnetic field $B$. In our study, the total Hamiltonian for a single particle is $\mathcal{H} =	\mathcal{H}_0	+\mathcal{H}^K$ where $\mathcal{H}_0$  denotes the kinetic and the Zeeman Hamiltonian and $\mathcal{H}^K$ the Rashba Hamiltonian.  $\mathcal{H}_0$ is given by
	\begin{align}
		\mathcal{H}_0 = \frac{p^2}{2m} +U(x,t)  - \hbar B (x,t)\cdot \zs\;, \label{tot_Ham}
	\end{align}
	where $p=-i\hbar \nabla$ is the momentum operator, $m$ the mass, and $x\in\mathbb{R}_x^3$ denotes the physical coordinate space. The field $B:\mathbb{R}_x^3\times\mathbb{R}_t\rightarrow \mathbb{R}^3$ represents some external time-dependent magnetic field. 
	The Rashba interaction is given by \cite{Manchon_15}
	\begin{align}
		\mathcal{H}^{K} = \hbar\left( p \wedge K\right) \cdot \zs \;,
		\label{Rash_Ham}
	\end{align}
	where $K:\mathbb{R}_x^3\times\mathbb{R}_t\rightarrow \mathbb{R}^3$ is a spatially dependent field describing the intensity and orientation of the Rashba interaction, 
	$\zs_i$ with $i=1,\ldots,3$ denote the Pauli matrices,
	and $\zs_0$ the two-by-two identity matrix. 
	We make the following technical assumptions concerning the fields 
	\begin{itemize}
		\item[AS:]  \emph{Assumption on the potential energy and Zeeman-Rashba interactions.} 
		We assume that the external magnetic and Rashba fields have bounded derivatives, in particular, $B_i \in C \left([0,T], W^{1,\infty}(\mathbb{R}^3_x)\right)$, and $K_i \in C \left([0,T], W^{2,\infty}(\mathbb{R}^3_x)\right)$, with $i=1,2,3$. Concerning the potential energy, we assume $U\in C^\infty_b \left(\mathbb{R}^3_x\times \mathbb{R}_t \right)$, where $C^\infty_b$ denotes the space of infinitely uniformly bounded differentiable functions.
	\end{itemize}
	We remark that Eq. \eqref{Rash_Ham} represents only the formal definition of the Rashba interaction, considered as a classical operator. A well-known fact concerning classical interactions expressed by products of the position and the momentum coordinates is the lack of uniqueness of well-defined quantum operators, whose classical limit agrees with Eq. \eqref{Rash_Ham}. Such a loss of one-to-one correspondence between classical phase-space functions and quantum operators arises from the non-commutativity between quantum position and momentum and is referred to as the quantization problem (for an analogous discussion in the case of wave propagating in nonuniform media see \cite{Morandi_25_KRM,Morandi_24_JPA}). In our work, we assume the Weyl quantization procedure, which associates any classical phase-space functions with the unique quantum operator that guarantees symmetric ordering of position and momentum operators. In this context, classical phase-space functions that represent physical observables are denoted symbols.
	%
	%
	The quantum dynamics of the particle is described by the Schr\"odinger problem
	$
	i\hbar \dpp{\zy}{t} = \textrm{Op}(\mathcal{H})\zy$, 
where $\textrm{Op}(\mathcal{H})$ denotes the Weyl operator that corresponds to the Hamiltonian symbol. The formal definition of the Weyl quantization is given by
\begin{align*}
	\left(\textrm{Op}( \mathcal{H}) \zy \right)(x)\doteq  \frac{1}{(2\pi\hbar)^{3}} \int_{\mathbb{R}^3_y\times\mathbb{R}^3_y} \mathcal{H}\left(\frac{x+y}{2},p \right)\zy (y ) e^{\frac{i}{\hbar}(x-y)p}\dif y\dif p\;,
\end{align*}
and the natural Hilbert space associated with the Schr\"odinger description of the quantum particle is $ L^2\left(\mathbb{R}^3;\mathbb{C}^2 \right)    $. 
\begin{lemma}\label{Lemma_ess_sa_H}
	Under AS, 
	the operator $ \textrm{Op}\left(\frac{p^2}{2m} +U(x) +\hbar \left(p\wedge K -B  \right) \cdot \zs  \right) $ is essentially self adjoint on
	$C_0^\infty(\mathbb{R}^3,\mathbb{C}^2)$.
\end{lemma}
\proof
Straightforward calculations show that
\begin{align*}
	\textrm{Op}\left(\frac{p^2}{2m} +U +\hbar \left(p\wedge K -B\right) \cdot \zs  \right) \zy    =& \left(-\frac{\hbar^2 \Delta  }{2m} + U  -\hbar B \cdot \zs \right)\zy \\
	&-i\hbar \left( \zs \cdot\left(\nabla \zy\wedge K \right)  +\frac{1}{2}\left(\nabla\wedge K \right) \cdot \zs  \zy \right).
\end{align*}
It is easy to verify that the operators $-\frac{\hbar^2 \Delta  }{2m} + U  -\hbar B \cdot \zs -i\hbar \left( \zs \cdot\left(\nabla \zy\wedge K \right)  +\frac{1}{2}\left(\nabla\wedge K \right) \cdot \zs  \zy \right)$ are symmetric in $C_0^\infty(\mathbb{R}^3)^3$ with respect to the Hilbert product $\langle\zy,\zf \rangle \doteq \sum_{r=1}^2\int_{\mathbb{R}^3}\overline{\zy}_r(x)\zf_r(x)\dif x$. Under the hypothesis on the electric and magnetic potentials and on the Rashba field, standard estimates ensure that the operator $ U  -\hbar B \cdot \zs -  \frac{i\hbar}{2}\left(\nabla\wedge K \right) \cdot \zs  $ is bounded with respect to the Laplacian operator, with a relative bound smaller than one. The operator $  -i\hbar  \zs \cdot\left(\nabla \zy\wedge K \right) $ can be estimated as follows
\begin{align*}
	\left\| -i\hbar  \zs \cdot\left(\nabla \zy\wedge K \right)\right\|^2_{L^2(\mathbb{R}^3,\mathbb{C}^2)} \leq& \hbar \max_{i} \| K^2_i \|_{L^\infty(\mathbb{R}^3)}  \sum_i \left\| \dpp{\zy}{x_i} \right\|_{L^2(\mathbb{R}^3,\mathbb{C}^2)}^2 \\ \leq &   \hbar \max_{i} \| K^2_i \|_{L^\infty(\mathbb{R}^3)}  \left\| (b+a|\zh|^2) \widetilde{\zy} \right\|_{L^2(\mathbb{R}^3,\mathbb{C}^2)}^2\\
	\leq &   \hbar \max_{i} \| K^2_i \|_{L^\infty(\mathbb{R}^3)} \left( b\left\|  \zy \right\|_{L^2(\mathbb{R}^3,\mathbb{C}^2)}^2+ a\left\| \Delta \zy \right\|_{L^2(\mathbb{R}^3,\mathbb{C}^2)}^2\right).
\end{align*}
where $\widetilde{\zy} \doteq \mathcal{F}_{x\rightarrow\zh}$ denotes the Fourier transform and we have used that for any $a>0$ arbitrarily small, there exists $b>0$ such that $|\zh|\leq b+a|\zh|^2$, for any $\zh\in \mathbb{R}^3$. The conclusion follows from the Kato-Rellich theorem. \endproof
\noindent 
In the remaining part of the study, we will ambient the analysis of the optimal control problem of the quantum dynamics in the Wigner-Moyal formalism, which refers to the general framework of a statistical mixture of quantum particles.
%
%
The Wigner pseudo-distribution function associated to the density matrix $\widehat{\zr}  =  |\zy \rangle \langle \zy|$, is a two-by-two matrix whose elements are defined as follows
\begin{align}
	f_{ij}\left(x,p\right)\doteq &\frac{1}{(2\pi\hbar)^{3}} \textrm{Op}^{-1}(\widehat{\zr})
	=  \frac{1}{(2\pi)^{3}}	\int_{\mathbb{R}^{3}} \zr_{ij} \left(x+\frac{\hbar{ \zh}}{2} ,x-\frac{\hbar{ \zh}}{2}\right)  e^{-ip \zh}\dif { \zh}\;,\label{def_W_op}
\end{align}
where $\zr_{ij}(x,x') = \zy_i(x)\overline{\zy}_j(x')$ and it is assumed $\zy\in L^2\left(\mathbb{R}^3,\mathbb{C}^2 \right)    $.
We denote by $f=\mathcal{W}^{\hbar}[\zy]$ the operator that associates the wave function $\zy$ with the corresponding hermitian matrix of Wigner functions, with elements given by Eq. \eqref{def_W_op}. 
We define the Hilbert space of phase-space functions in L2 taking values in the space of two-by-two complex matrices
\begin{align*}
	\mathrm{HS} \left(\mathbb{R}^3_x\times \mathbb{R}^3_p\right)=\left\{ f : \mathbb{R}^3_x \times \mathbb{R}^3_p \rightarrow \mathbb{C}^4 : \mathrm{tr} \int_{\mathbb{R}_x^3\times \mathbb{R}^3_p}  f^\dag f\dif x\dif p<\infty \right\}\;,
\end{align*}
where $\dag$ denotes the hermitian conjugate and $\mathrm{tr} $ the trace. The associated Hilbert inner product is given by $\left\langle f,h \right\rangle \doteq  \mathrm{tr} \int_{\mathbb{R}_x^3\times \mathbb{R}^3_p}  f^\dag h \dif x\mathrm{d} p$. In the mathematical analysis of the OC problem, we will consider the H1 Sobolev space with width, defined by the following norm
\begin{align*}
	\| f \|_{H^1_p} \doteq
	& \left(\sum_{i=0}^3 \int_{\mathbb{R}_x^3\times \mathbb{R}^3_p} \left[\left(1+ |p|^2 \right) f_i^2 +\sum_{j=1}^3\left( \left|\dpp{f_i}{x_j} \right|^2 + \left|\dpp{f_i}{p_j} \right|^2  \right) \right]\dif x\dif p\right)^{\frac{1}{2}}\;.
\end{align*}
We recall few elementary facts concerning the Wigner matrix function. The L2 norm of the Wigner matrix function scales as $\hbar^{-\frac{3}{2}}$ with respect to the squared L2 norm of the wave function, $
\left\| f \right\|_{\textrm{HS}\left(\mathbb{R}^3_x\times \mathbb{R}^3_p\right)}=\frac{1}{\hbar^{3/2}} \left\| \zy \right\|_{L^2\left(\mathbb{R}^3_x,\mathbb{C}^2\right)}^2$.
Moreover, the zeroth-order moment of the Wigner matrix function provides the squared L2 norm of the associated wave function, $\left\| \zy \right\|_{L^2\left(\mathbb{R}^3_{x};\mathbb{C}^2\right)}^2
=\textrm{tr} \int_{\mathbb{R}^6_{x,p}}  f \left(x,p\right) \dif x \dif p$, as expected. 
%
The evolution equation for the Wigner matrix is obtained by applying the inverse Weyl map $\textrm{Op}^{-1}$ to the von Neumann equation
$	i\hbar \dpp{\widehat{\zr}}{t}  = \left[\widehat{\mathcal{H}},\widehat{\zr} \right] $. We obtain
\begin{align}
	i\hbar \dpp{f}{t} =& \left[ \left( \frac{p^2}{2m} +U +\hbar\left( p \wedge K -  B\right) \cdot \zs \right),f \right]_{\#} \;,\label{Wig_dyn_Moy}
\end{align}
where $\#$ denotes the Moyal product, whose explicit form is given in Eq. \eqref{Moyal_prod}.
The Moyal product between two phase-space functions assumes a relatively simple expression if one of the two symbols depends only on the spatial or momentum coordinate. For instance, in the case of the potential energy $U=U(x)$, which depends only on the position variables, the Moyal commutators and anti-commutators can be written as $\left[ U ,F \right]_{\#}  \doteq U\#F-F\#U   =i\hbar \Theta^-_{U}[F] $ and $\left\{ U ,F \right\}_{\#} \doteq U\#F+F\#U = \Theta^+_{U}[F]$, respectively. We have defined 
\begin{align}
	\Theta^\pm_{U}[F](x,p)  \doteq&    \frac{ 1}{(2\pi)^{3}}  \int_{\mathbb{R}^{6}} \delta_\pm U \left(x ,\zh  \right)F \left(x,p' \right) e^{- i\zh(p -p') }    \dif \zh   \dif p'\;, \label{Theta_pm}
\end{align}
with
\begin{align}
	\delta_- U \left(x ,\zh  \right) \doteq& \frac{1}{i\hbar}\left[U \left(x +\frac{\zh\hbar}{2}  \right)-U \left(x -\frac{\zh\hbar}{2}  \right)\right] \label{delta_k_-}\\
	\delta_+ U \left(x ,\zh  \right) \doteq & U \left(x +\frac{\zh\hbar}{2}  \right)+U \left(x -\frac{\zh\hbar}{2}  \right)\label{delta_k_+}\;.
\end{align}
We remark that if $U\in L^\infty(\mathbb{R}^3_x)$, $\Theta^\pm_U$ operators are bounded in $\textrm{HS}\left(\mathbb{R}^3_x\times \mathbb{R}^3_p\right)$. In fact,
\begin{align*}
	\left\|\Theta^-_{U}[F]\right\|_{\textrm{HS}\left(\mathbb{R}^3_x\times \mathbb{R}^3_p\right)}   =&  \left\|\mathcal{F}_{p\rightarrow \zh}^{-1}\Theta^-_{U}[F]\right\|_{\textrm{HS}\left(\mathbb{R}^3_x\times \mathbb{R}^3_\zh\right)} = \frac{1}{\hbar}\left\|\left[U \left(x +\frac{\hbar \zh}{2}  \right)-U \left(x -\frac{\hbar \zh}{2}  \right) \right] \mathcal{F}_{p\rightarrow \zh}^{-1}F\right\|_{\textrm{HS}\left(\mathbb{R}^d_x\times \mathbb{R}^d_\zh\right)} \\
	\leq&    \frac{2}{\hbar}\left\|U\right\|_{L^\infty\left(\mathbb{R}^3_x\right)} \left\|F\right\|_{\textrm{HS}\left(\mathbb{R}^3_x\times \mathbb{R}^3_\zh\right)} ,
\end{align*}
and similar for the $\Theta^+$ operator.

In order to study Eq. \eqref{Wig_dyn_Moy}, it is convenient to decompose the Wigner matrix on the Pauli basis. We write $f=\sum_{i=0}^3 f_i \zs_i$, where the projection coefficients are given by $f_i=\frac{1}{2}\mathrm{tr} \left\{\zs_i f \right\}$. In particular, if $f$ is a hermitian matrix, the coefficients $f_i$ are real. The results of our mathematical analysis concerning Eq.  \eqref{Wig_dyn_Moy} are presented in the following
\begin{theorem}\label{teor_wellpos_W_dyn}
	Let $T>0$, and assume AS. 
	The Wigner transform $f\doteq \mathcal{W}^{\hbar}[\zy]$ of the solution of the Schr\"odinger equation with Hamiltonian given by Eq. \eqref{tot_Ham} and initial condition $\zy_0\in L^2(\mathbb{R}^3,\mathbb{C}^2)$, and $\mathcal{W}^{\hbar}[\zy_0]\in H^1_p$,  is the unique strong solution of Wigner problem
	\begin{align}
		\left\{\begin{array}{lll}
			\disp	\dpp{f_0}{t} + \frac{p}{m} \cdot\nabla_xf_0 -\Theta^-_U [f_0]=&  \disp \sum_{i=1}^3\left(A^+[f_i]_i  - \hbar \Theta^-_{ B_i}[f_i]\right) &\textrm{on } \mathbb{R}_x^3\times\mathbb{R}^3_p\times[0 ,T] \\[6pt]
			\disp	\dpp{f_k}{t} + \frac{p}{m} \cdot\nabla_xf_k -\Theta^-_U [f_k]=&\disp A^+[f_0]_k   + \sum_{i,j=1}^3 \ze_{ijk}  \, A^-[f_j]_i & \\[6pt]
			&   \disp   -   \hbar\Theta^-_{ B_k}[f_0]-  \sum_{i,j=1}^3 \ze_{ijk}\, \Theta^+_{ B_i}[f_j] &\textrm{on } \mathbb{R}_x^3\times\mathbb{R}_p^3\times[0 ,T] \\[6pt]
			\disp   \left.f \right|_{t=0} =\mathcal{W}^\hbar [\zy_0] &&\textrm{on } \mathbb{R}_x^3\times\mathbb{R}_p^3
		\end{array}
		\right.\label{Wig_f0_fk_01}
	\end{align}
	the $k$ index in the second line takes the values $k=1,2,3$, the $A^\pm$ operators are  defined as
	\begin{align}
		A^+[h]_k=	& \hbar \sum_{i,j=1}^3  \ze_{ijk}\left(   p_i 	\Theta^-_{K_j}[h]   -  \frac{1}{2}  \Theta^+_{K_j}\left[\dpp{h}{x_i}\right] \right)\\
		A^-[h]_k=	& \sum_{i,j=1}^3  \ze_{ijk}\left(    p_i 	\Theta^+_{K_j}[h]   +  \frac{\hbar^2}{2}  \Theta^-_{K_j}\left[\dpp{h}{x_i}\right] \right)\;,
	\end{align}
	and the pseudo-differential operators $\Theta^\pm$ are defined in Eq. \eqref{Theta_pm}.
	Moreover, $f\in C\left([0,T], H^1_p\right)$, and for all $t\in[0,T] $,  $\|f(.,.,t)\|_{L^2\left( \mathbb{R}^3_x\times\mathbb{R}^3_p,\mathbb{R}^4 \right)}=\|f(.,.,0)\|_{L^2\left( \mathbb{R}^3_x\times\mathbb{R}^3_p,\mathbb{R}^4 \right)}$ and $f $ satisfies 
	\begin{align}
		\| f(\cdot,\cdot,t) \|_{H^1_p} \leq&  \| f(\cdot,\cdot,0) \|_{H^1_p} e^{C \int_0^t \left(\| K(.,t)\|_{W^{2,\infty}(\mathbb{R}_x^3)^3} +\| U(.,t)\|_{W^{1,\infty}(\mathbb{R}_x^3)}+\| B(.,t)\|_{W^{1,\infty}(\mathbb{R}_x^3)^3} \right)\dif t} \;.\label{H1_bound_sol_wig}
	\end{align}
\end{theorem}
Before providing the proof, we discuss some formal limits and particular cases of the Wigner equations. Due to the complexity of the full quantum Wigner dynamics, in physical applications, it may be convenient to approximate the full quantum dynamics with the corresponding classical limit, obtained formally by expanding the operators with respect to $\hbar$, and letting $\hbar $ go to zero. By Taylor expansion, the pseudo-differential operators can be approximated as
$	\Theta^-_{h}[f] \stackrel{\hbar\rightarrow 0}{\longrightarrow }\nabla_x h \nabla_p f +o(\hbar)$ and $
\Theta^+_{h}[f]  \stackrel{\hbar\rightarrow 0}{\longrightarrow } 2 h  f+o(\hbar)$. Equation  \eqref{Wig_f0_fk_01} simplifies 
\begin{align}
	\dpp{f_0}{t}  +  \frac{p}{m} \cdot\nabla_xf_0 -\nabla_xU \cdot\nabla_p f_0 =&\hbar\left(\sum_{i,j,k=1}^3  \ze_{ijk}\left( p_i 	\nabla_x K_j\cdot \nabla_p f_k     \right) +  K\cdot\left(   \nabla_x \wedge f  \right)    -  \sum_{i=1}^3 \nabla_x B_i \cdot\nabla_p f_i\right) \label{form_hbar_exp_f0} \\
	\dpp{f_k}{t} + \frac{p}{m} \cdot\nabla_xf_k -\nabla_xU \cdot\nabla_p f_k =& 2    \left[( p\wedge K - B )\wedge f    \right]_{k} \nn \\
	&+\hbar \left(\sum_{i,j=1}^3  \ze_{ijk}\left(     p_i 	\nabla_x K_j\cdot \nabla_p f_0      \right)+   \left(   K\wedge \nabla_x   \right)_k f_0   - \nabla_x B_k\cdot \nabla_p f_0\right)\nn \\
	& + \frac{ \hbar^2}{2}  \sum_{i,j,l,s=1}^3 \ze_{ils} \;    \ze_{ijk}  \nabla_x K_s\cdot \nabla_p \left( \dpp{f_j}{x_l}  \right),& \label{form_hbar_exp_fk}
\end{align}
where $k=1,2,3$. The correspondence between the full quantum Wigner dynamics and the Liouville-Rashba equation in the classical limit will be discussed in Sec. \ref{Sec_cl_lim}. 
Typically, realistic experimental setups investigating the quantum dynamics of electrons in semiconductors with Rashba effect, consider the presence of uniform external magnetic field, and employ materials whose Rashba coefficient $K$ is uniform in space. In such cases, without any further simplifying assumptions, the quantum evolution equations simplify  
\begin{align}
	&\disp	\dpp{f_0}{t} + \frac{p}{m} \cdot\nabla_xf_0 -\Theta^-_U [f_0]=    \hbar  K\cdot\left(   \nabla_x \wedge f  \right)  &\textrm{on } \mathbb{R}_x^d\times\mathbb{R}^d_p\times[0 ,T] \\[6pt]
	&\disp	\dpp{\vec{f}}{t} + \frac{p}{m} \cdot\nabla_x \vec{f} -\Theta^-_U [\vec{f}]=  2    ( p\wedge K - B )\wedge \vec{f} +\hbar \left(   K\wedge \nabla_x   \right) f_0   &\textrm{on } \mathbb{R}_x^d\times\mathbb{R}_p^d\times[0 ,T] 
	\label{Wig_KB_uni_01}
\end{align}
where, to compact notations, we have indicated by a vector the three cartesian components of the Wigner matrix ${\vec{f}}_i=\frac{1}{2} \mathrm{tr}\left\{ f \zs_i \right\} $ with $i=1,2,3$. We now move on to the proof of Th. \ref{teor_wellpos_W_dyn}.
\proof \emph{of Th. \ref{teor_wellpos_W_dyn}.} 
The Wigner equation \eqref{Wig_f0_fk_01} is obtained by applying the definition of Moyal product in Eq. \eqref{Wig_dyn_Moy}. The details of the derivation are given in Appendix \ref{App_der_eq_wig}. According to Lemma \ref{Lemma_ess_sa_H}, the Hamiltonian operator $\textrm{Op}(\mathcal{H})$ generates the unitary group $e^{-\frac{i}{\hbar}\textrm{Op}(\mathcal{H})t} $ in $L^2(\mathbb{R}^3_x,\mathbb{C}^2)$. The Wigner transform of $\zy(t)=e^{-\frac{i}{\hbar}\textrm{Op}(\mathcal{H})t}\zy_0$ provides the unique solution of Eq.  \eqref{Wig_dyn_Moy} or, equivalently, of Eq.\eqref{Wig_f0_fk_01}.
The conservation of the L2 norm of the solution is a direct consequence of the equality $
\left\| \mathcal{W}^\hbar [\zy] \right\|_{\textrm{HS}\left(\mathbb{R}^d_x\times \mathbb{R}^d_p\right)}=\frac{1}{\hbar^{3/2}} \left\| \zy \right\|_{L^2\left(\mathbb{R}^d_x,\mathbb{C}^2\right)}^2$, mentioned before. 

The derivation of Eq. \eqref{H1_bound_sol_wig} requires cumbersome calculations. Here, for the sake of simplicity, we limit ourselves to considering the Rashba Hamiltonian $\mathcal{H}^K$. Consequently, in the rest of the proof, we refer to the simplified equation $i\hbar\dpp{f}{t}=\left[\mathcal{H}^K,f \right]_\#$. The remaining terms can be treated in the same way. The reader interested in more details may refer to \cite{Morandi_SIAM_25} in which the optimal control problem of the quantum dynamics in the absence of magnetic and Rashba interaction is considered. The proofs can be easily adapted to derive the bounds associated to the kinetic energy, electric and magnetic fields.

We define the linear operator $ \Omega \doteq   1+ \sum_{i=1}^3 \left( p_{i}^2 - \dsp{}{x_{i}} -\dsp{}{p_{i}} \right)$. From the definition of the scalar product $
\left\langle 	f   ,h \right\rangle \doteq \frac{1}{2}\mathrm{tr}  \int_{\mathbb{R}^{6}} f^\dag h  \dif p\dif x  $,  we have $\|f \|^2_{H^1_p} =\left\langle 	f   ,\Omega f \right\rangle $, where we used $f=f^\dag$. We have
\begin{align}
	\dpt{}{t}\left\langle f, \Omega f \right\rangle =& \frac{i}{\hbar}    \left\langle  \left(\left[\mathcal{H}^K,\Omega f  \right]_\#-\Omega \left[\mathcal{H}^K,f \right]_\#\right)  ,f \right\rangle \nn \\ =&	   \sum_{i,j,k=1}^3  \ze_{ijk} \left(	\mathcal{A}^{ij}_{\zO}(f_0,f_k) +	\mathcal{A}^{ij}_{\zO}(f_k,f_0) \right)+\sum_{i,j,k=1}^3 \ze_{ijk} \sum_{r,s=1}^3\ze_{ksr}	\mathcal{B}^{ij}_{\zO}(f_s,f_r)\;,\label{dt_fOf}
\end{align}
where we have defined
\begin{align*}
	\mathcal{A}^{ij}_{\Omega}(f_0,f_k) \doteq&  	i   \int_{\mathbb{R}^{6}}	\left( \left[ p_i K_j ,\Omega f_0 \right]_\# - 	\Omega  \left[ p_i K_j ,f_0 \right]_\#\right) f_k   \dif p\dif x \\
	\mathcal{B}^{ij}_{\Omega}(f_s,f_r) \doteq&  	- 	 	 \int_{\mathbb{R}^{6}} \left(\left\{ p_i K_j  , \Omega f_s \right\}_\#  -\Omega   \left\{ p_i K_j  ,  f_s \right\}_\#  \right) f_r \dif p\dif x   \;.
\end{align*}
					We define $\Omega_0^i\doteq p^2_i$, $\Omega_1^i\doteq -\dsp{}{p_i}$,  $\Omega_2^i\doteq-\dsp{}{x_i}$ and $\Omega=1+\sum_{i=1}^3 \left(\Omega_0^i+\Omega_1^i+\Omega_2^i\right)$. We analyze each term separately. We start with $\Omega_0^i$. Using Eq. \eqref{moy_prod_pKf_01}, we have
					\begin{align}
						\mathcal{A}^{ij}_{p^2_\mu}(f_0,f_k) \doteq&   	  \hbar \int_{\mathbb{R}^{6}}	\left( 	- p_i  \Theta^-_{K_j}[p^2_\mu f_0]
						+ p_i p^2_\mu \Theta^-_{K_j}[f_0]   - \frac{1}{2} p^2_\mu \Theta^+_{K_j}\left[\dpp{f_0}{x_i}\right]    + \frac{1}{2}  \Theta^+_{K_j}\left[p^2_\mu \dpp{f_0}{x_i}\right]  \right) f_k   \dif p\dif x    \:. \label{A_p_mu}
					\end{align}
					%
						We estimate the first two terms in Eq. \eqref{A_p_mu}.  Using $\left({p_\mu'}^2- p_\mu^2\right)e^{- i\zh(p -p') } =-2i p_\mu \dpp{}{\zh_\mu} e^{- i\zh(p -p') } - \dsp{}{\zh_\mu} e^{- i\zh(p -p') } $ and integrating by parts, we have
						\begin{align*}
							\Theta^-_{K}[p^2_\mu f_0]    -	p^2_\mu \Theta^-_{K}[ f_0]
							=& \frac{ 1}{(2\pi)^{3}}	\int_{\mathbb{R}^{6}}      \delta_- K \left(x, \zh \right)  \left(-2i p_\mu \dpp{}{\zh_\mu} e^{- i\zh(p -p') } - \dsp{}{\zh_\mu} e^{- i\zh(p -p') } \right)  f \left(x,p' \right)     \dif \zh   \dif p'  \\
							=& 	p_\mu \Theta^+_{\dpp{K}{x_\mu} }[ f_0]    -	 	\frac{\hbar^2}{4}\Theta^-_{\dpp{K}{x_\mu}}[ f_0]\;,
						\end{align*}
						where the definitions of $\delta_- K$ and $\delta_+ K$
						are given in Eq. \eqref{delta_k_-}.
							Cauchy-Schwartz inequality gives
							\begin{align*}
								&	  \left| \int_{\mathbb{R}^{6}} p_i 	 p_\mu \Theta^+_{\dpp{K}{x_\mu} }[ f_0]  f_k   \dif p\dif x \right|\leq   \left\|  p_\mu \Theta^+_{\dpp{K}{x_\mu} }[ f_0] \right\|_{L^2 (\mathbb{R}^{6}_{xp})} \left\| p_i 	 f_k    \right\|_{L^2 (\mathbb{R}^{6}_{xp})}
							\end{align*}
							and
							\begin{align*}
								\left\|  p_\mu \Theta^+_{\dpp{K}{x_\mu} }[ f_0] \right\|_{L^2 (\mathbb{R}^{6}_{xp})} 	=&	  \left\|  \dpp{}{\zh_\mu} \delta^+_{\dpp{K}{x_\mu} }  \widehat{f_0} \right\|_{L^2 (\mathbb{R}^{6}_{x\zh})}
								\\
								\leq &	    \left(\hbar \left\| \dsp{K}{x_\mu}  \right\|_{L^\infty (\mathbb{R}^{d}_{x})^3} \left\|   f_0 \right\|_{L^2 (\mathbb{R}^{6}_{xp})} + 2 \left\| \dpp{K}{x_\mu}  \right\|_{L^\infty (\mathbb{R}^{3}_{x})^3} \left\|    p_\mu f_0 \right\|_{L^2 (\mathbb{R}^{6}_{xp})} \right)\;,
							\end{align*}
							where $\widehat{f_0} = \mathcal{F}_{p\rightarrow \zh} (f_0) (x,\zh)$ denotes the Fourier transform with respect to the momentum variable.
							Similarly,
							\begin{align*}
								&\frac{\hbar^3}{4} \left| \int_{\mathbb{R}^{6}} p_i  	 	\Theta^-_{\dpp{K}{x_\mu}}[ f_0]   f_k   \dif p\dif x \right| \leq \frac{\hbar^3}{4} \left| \int_{\mathbb{R}^{6}} \Theta^-_{\dpp{K}{x_\mu} }[ f_0] p_i 	 f_k   \dif p\dif x \right|\leq \frac{\hbar^3}{4}  \left\|    \Theta^-_{\dpp{K}{x_\mu} }[ f_0] \right\|_{L^2 (\mathbb{R}^{6}_{xp})} \left\| p_i 	 f_k    \right\|_{L^2 (\mathbb{R}^{6}_{xp})} \\
								=&	\frac{\hbar^3}{4} \left\|  \delta^-_{\dpp{K}{x_\mu} }  \widehat{f_0} \right\|_{L^2 (\mathbb{R}^{6}_{x\zh})} \left\| p_i 	 f_k    \right\|_{L^2 (\mathbb{R}^{6}_{xp})}	\leq 	 	\frac{\hbar^2}{2} \left\| p_i 	 f_k    \right\|_{L^2 (\mathbb{R}^{6}_{xp})}  \left\| \dpp{K}{x_\mu}  \right\|_{L^\infty (\mathbb{R}^{d}_{x})^3} \left\|   f_0 \right\|_{L^2 (\mathbb{R}^{6}_{xp})}\;.
							\end{align*}
							The remaining terms of Eq. \eqref{A_p_mu} can be treated in a similar way. We write $	p^2_\mu \Theta^+_{K_j}\left[\dpp{f_0}{x_i}\right]    - \Theta^+_{K_j}\left[p^2_\mu \dpp{ f_0}{x_i}\right]   =   \hbar^2	p_\mu \Theta^-_{\dpp{K}{x_\mu} }\left[  \dpp{ f_0}{x_i} \right]    +	\frac{\hbar^2}{4}\Theta^+_{\dsp{K}{x_\mu}} \left[  \dpp{ f_0}{x_i} \right]  $, and we have the estimates
							\begin{align*}
								\left| \int_{\mathbb{R}^{6}} \hbar^2	p_\mu \Theta^-_{\dpp{K}{x_\mu} }\left[  \dpp{ f_0}{x_i} \right]     f_k   \dif p\dif x \right|
								%
								\leq &	 2	 \hbar  \left\| \dpp{K}{x_\mu}  \right\|_{L^\infty (\mathbb{R}^{3}_{x})^3} \left\|     \dpp{ f_0}{x_i}     \right\|_{L^2 (\mathbb{R}^{6}_{xp})} \left\| p_\mu 	 f_k    \right\|_{L^2 (\mathbb{R}^{6}_{xp})}  \\
								\left| \int_{\mathbb{R}^{6}}	\frac{\hbar^2}{4}	  \Theta^+_{\dsp{K}{x_\mu} }\left[  \dpp{ f_0}{x_i} \right]     f_k   \dif p\dif x \right|
								%
								\leq &	 \frac{\hbar^2}{2}	\left\| \dsp{K}{x_\mu}  \right\|_{L^\infty (\mathbb{R}^{3}_{x})^3} \left\|     \dpp{ f_0}{x_i}     \right\|_{L^2 (\mathbb{R}^{6}_{xp})}  \left\|  	 f_k    \right\|_{L^2 (\mathbb{R}^{6}_{xp})}\;.
							\end{align*}
							In conclusion, we have obtained 
							\begin{align*}
								\left|	\mathcal{A}^{ij}_{p^2_\mu}(f_0,f_k)\right| 	\leq &\left(  \left\| \dsp{K}{x_\mu}  \right\|_{L^\infty (\mathbb{R}^{3}_{x})^3}  +2  \left\| \dpp{K}{x_\mu}  \right\|_{L^\infty (\mathbb{R}^{3}_{x})^3} \right)\\ & \left(\hbar	 \left\| p_i 	 f_k    \right\|_{L^2 (\mathbb{R}^{6}_{xp})}  + \hbar^2	 \left\|  	 f_k    \right\|_{L^2 (\mathbb{R}^{6}_{xp})}  \right) \left( \hbar \left\|   f_0 \right\|_{L^2 (\mathbb{R}^{6}_{xp})} +  \left\|    p_\mu f_0 \right\|_{L^2 (\mathbb{R}^{6}_{xp})} +  \left\|     \dpp{ f_0}{x_i}     \right\|_{L^2 (\mathbb{R}^{6}_{xp})} \right)\\
								\leq &C  \left\| K\right\|_{W^{2,\infty} (\mathbb{R}^{d}_{x})^3}   \|f\|_{H_p^1}^2\;.
							\end{align*}
							The term $\mathcal{B}$ in Eq. \eqref{dt_fOf} containing the Moyal anti-commutators can be treated in the same way, and we obtain
							\begin{align*}
								\dpt{}{t}	\left\langle f(t), \sum_{\mu } p^2_\mu f (t) \right\rangle  
								\leq C  \left\| K\right\|_{W^{2,\infty} (\mathbb{R}^{3}_{x})^3}   \|f\|_{H_p^1}^2\;.
							\end{align*}
							Concerning $\Omega_1^\mu=-\dsp{}{p_\mu}$, we observe that the following simplifications hold
							\begin{align*}
								&\left[  p_i K_j , \dsp{}{p_\mu} f_0 \right]_{\#}	-\dsp{}{p_\mu}	\left[  p_i K_j , f_0 \right]_{\#} =  	-i 2 \hbar\;    \zd_{i,\mu}\Theta^-_{K_j}\left[\dpp{f_0}{p_\mu} \right]  \\
								&\left\{ p_i K_j , \dsp{}{p_\mu} f_0 \right\}_{\#}	-\dsp{}{p_\mu}	\left\{  p_i K_j , f_0 \right\}_{\#} =  	-2 \zd_{i,\mu}   \Theta^+_{K_j}\left[\dpp{f_0}{p_\mu} \right]\;,
							\end{align*}
							where $\zd$ denotes the Kronecker's delta.
							Proceeding as before, we obtain
							\begin{align*}
								&\dpt{}{t}\left\langle f,  -\sum_{\mu} \dsp{f}{p_\mu}   \right\rangle  \leq C   \max_j    \left\| K_j  \right\|_{L^\infty (\mathbb{R}^{3}_{x})^3}
								\|f\|_{H_p^1}^2\;.
							\end{align*}
							Finally, for $\Omega^\mu_2=-\dsp{}{x_\mu}$ we have $	\dpt{}{t}\left\langle f,  -\sum_{\mu}\dsp{f}{x_\mu}   \right\rangle  \leq C   \max_\mu  \left(  \left\| \dsp{K}{x_\mu}  \right\|_{L^\infty (\mathbb{R}^{3}_{x})^3}  + \left\| \dpp{K}{x_\mu}  \right\|_{L^\infty (\mathbb{R}^{3}_{x})^3}  \right)  \left\|    f   \right\|_{H^1_p } ^2$. Details are given in Appendix \ref{App_est_d_xmu}. Putting together, we have obtained $\dpt{}{t}	 \|f\|_{H_p^1}^2
							\leq C  \left\| K\right\|_{W^{2,\infty} (\mathbb{R}^{3}_{x})^3}   \|f\|_{H_p^1}^2$.
							As we have indicated before, in our estimations we have neglected the Hamiltonian terms $\mathcal{H}' = \frac{p^2}{2m} +U(x)  - \hbar B (x)\cdot \zs$.  Including such terms is straightforward and leads to slight modifications of the previous estimate.   Careful investigation gives 
							\begin{align*}
								\dpt{}{t}	 \|f\|_{H_p^1}^2
								\leq C \left( \left\| K\right\|_{W^{2,\infty} (\mathbb{R}^{3}_{x})^3} +\| U\|_{W^{1,\infty}(\mathbb{R}_x^3)}+\| B \|_{W^{1,\infty}(\mathbb{R}_x^3)^3} \right) \|f\|_{H_p^1}^2\;.
							\end{align*}
							Gr\"onwall's Lemma provides Eq. \eqref{H1_bound_sol_wig}.
							\endproof
							\section{Optimal Control of the quantum dynamics }
							
							In this section, we discuss the optimal control strategy that we apply to the quantum system. Our main goal is to dynamically control the quantum state of the system by using an external field whose shape may be modulated by few parameters. We introduce a set of $d$ time-dependent parameters represented by a function $u:\mathbb{R}_t\rightarrow \mathbb{R}^d$. As already anticipated in AS, for simplicity, we assume that the total potential energy acting on the system, which is the sum of some external potential representing the material plus the control, is a smooth function with bounded derivatives i. e. $U(x,u) \in C^\infty_b \left(\mathbb{R}^3\times \mathbb{R}^d \right)$, where $C_b^\infty $ denotes the $C^\infty $ space with bounded derivatives of all orders. Consistently with the previous assumptions on the total potential energy, we require that the temporal evolution of the potential energy originates from the control parameters. In various applications, ensamble of quantum particles may be controlled by trapping the system inside a harmonic trap. By designing the oscillator strength and the position of the minimum of the trap potential, it is possible to achieve some control of the quantum state of the system \cite{Morandi_SIAM_25,Cavaliere_19}. 
							Dipolar or harmonic controls have the technological advantage that they may be implemented with conventional electronics, so that the trapping potential can be modified by simply varying the gate voltage of a gate electrode.  
							For our purposes, it is not necessary to fix the explicit form of the controlled potential. The optimality conditions on the control parameters $u$, i. e. the set of equations that the time dependent control should satisfy in order to achieve optimal control of our system, are independent from the details of the parametrization of the control potential.  %
							In order to set up an optimal control procedure, it is necessary to specify the kind of control on the quantum state that we would like to achieve and, eventually, the cost on the control apparatus that we would like to minimize. First, we specify the goal of our control procedure. The aim of our strategy is to transfer a quantum particle from some initial position to a target position expressed by the classical phase-space coordinates $(x_T,p_T)$. The justification of the relevance of such a strategy in the case of the transport of cold atoms by optical tweezers, the discussion of the physical meaning of localizing a quantum state around a target point in the quantum phase-space and the link with the analogous classical transport, can be found in \cite{Morandi_PRA_25}. The spin of the particles constitutes an additional degree of freedom, which can also be controlled. The goal of transferring the particle at the final time $T$ with the highest possible probability, to the target position with the desired momentum and spin oriented along a fixed direction $d_T$, can be achieved by minimizing the following goal functional
							\begin{align}
								\Phi(f)=&\frac{1}{2} \mathrm{tr} \int f_{T}(x,p)f(x,p,T)\dif x \dif p.\label{goal_funct}
							\end{align}
							where the function $f_T$ measures the distance to the target and is chosen as follows 
							\begin{align}
								f_{T} (x,p)= \zs_0 \left[\frac{\nu_x}{2}|x -x_T|^2+\frac{\nu_p}{2} p^2  \right]- \nu_d \; d_T\cdot \zs , \label{target_funct}
							\end{align}
							where $(x_T,p_T)$ are the target coordinates in the phase-space and $d_T\in\mathbb{R}^3$ is the desired orientation of the dipole momentum at the final time. 
							We have introduced a set of positive widths $\nu_x,\nu_p,\nu_d$, which can be tailored to fix the importance to achieve some goal relative to the others.
							Typically, the control is achieved at the expense of some cost, which one would like to minimize and which physically represents the amount of external energy required by the control system. Furthermore, due to the technical implementation of the control system, there may be a technological or mechanical constraint that needs to be satisfied. In our model, we assume that the total cost associated with the control system can be described by the following functional    
							\begin{align}
								k(u)=&\frac{1}{2}\int_0^{T}\biggl[\gamma|u(t)|^2+ \gamma'\biggl|\dpt{u}{t}\biggr|^2\biggr]\dif t,\label{cost_func}
							\end{align}
							where $|\cdot |$ denotes the vector norm.
							Equation \eqref{cost_func} represents the electric energy of the control signals, and the widths $\zg,\zg'$ are positive quantities. From a mathematical point of view, $k(u)$ can be interpreted as a regularizing functional for control parameters. In fact, the employ of such a cost function in our OC procedures leads to controls belonging to H1.
							
							In synthesis, the OC problem is finalized to design a set of time-dependent parameters $u$ modulating the total electric potential applied to the system, with the aim of steering the evolution of a quantum particle and achieving two main goals: preparing some quantum state and minimizing costs. The main objective is to maximize the probability of finding the particle at the final time close to the position $x_T$, with momentum $p_T$ and spin direction aligned to the vector $d$. At the same time, the control should be implemented by using the minimum possible energy. From a mathematical point of view, the OC problem is formalized as follows
							\begin{align}
								\text{Q-OC: }\: 	\underset{ u }{\text{min}}\; \left[\Phi(f)+k(u)\right] &\quad  \textrm{s. t. }  \quad \textrm{Eq. } \eqref{Wig_f0_fk_01} \textrm{ holds true} \; .\label{OC_Q}
							\end{align}
							The OC problem can be formulated in terms of a variational problem, and the set of equations characterizing the stationary solutions is denoted an optimality condition. 
							A natural way to derive the optimality conditions of the OC problem of Eq. \eqref{OC_Q} is to apply a Lagranginn approach. We define the quantum Lagrangian 
							\begin{align*}
								\mathcal{L}\doteq&  \Phi + k+\frac{1}{2}\textrm{tr} \int_0^{T} \int_{\mathbb{R}^{6}} \left(	 \dpp{f}{t} - \frac{1}{i\hbar}\left[ \left( \frac{p^2}{2m} +U +\hbar \left( p \wedge K- B \right) \cdot \zs  \right),f \right]_{\#} \right)   h  \dif x\dif p \dif t\;.
							\end{align*}
							Reminiscent of the Lagrangian multipliers method for the stationary points associated to the problem of minimizing functionals subject to constraints, the quantum Lagrangian includes the evolution equation expressed in weak form. In the present context, the role of test functions is played by the so-called adjoint functions $h$ which generalize the definition of Lagrangian multipliers.
							The variations of the quantum Lagrangian with respect to all the dynamical variables (the Wigner function $f$, the adjoint function $h$ and the control $u$), provides the set of optimality conditions (see i. e. \cite{Morandi_PRA_25,Morandi_SIAM_25} for details). 
							The variation of the quantum Lagrangian with respect to the Wigner function $f$ provides the evolution equations for the adjoint variables. They solve the same equation as the Wigner pseudo-distribution $f$.
							\begin{align*}
								\left\{\begin{array}{ll}
									\disp i\hbar \dpp{h}{t} =  \left[ \left( \frac{p^2}{2m} +U  +\hbar\left( p \wedge K-   B \right) \cdot \zs  \right),h \right]_{\#} & \textrm{ on } \mathbb{R}_x^3\times\mathbb{R}_p^3\times[0 ,T
									]
									\\
									\left.h\right|_{t=T} =f_T  & \textrm{ on } \mathbb{R}_x^3\times\mathbb{R}_p^3
								\end{array}\right.
							\end{align*}
							We note an important loss of symmetry between the formulation of the Wigner and of the adjoint function problems. The initial datum of the Wigner function has a physical origin and represents the initial state of the quantum system. 
							On the other hand, the value of the adjoint function is prescribed in the final time $T$, the final condition depends only on the optimization strategy and coincides with the function $f_T$. 
							Taking the variation of the quantum Lagrangian with respect to the parameters $u$ we obtain the optimality condition for the control 
							\begin{align*}
								\zg' \dst{u_i}{t}-	\zg u_i   =&-   \frac{1}{2} \textrm{tr} \int_{\mathbb{R}_x^3\times\mathbb{R}_p^3}   \Theta^-_{\dpp{U}{u_i}} [h] \;f    \dif x  \dif p& i=1,\ldots,d.
							\end{align*}
							In summary, the set of non-linear equations associated to the stationary point of the OC problem of Eq. \eqref{OC_Q} is
							\begin{align} 
								\left\{ \begin{array}{ll}
									\disp	
									i\hbar \dpp{f}{t} =  \left[ \left( \frac{p^2}{2m} +U+\hbar\left( p \wedge K -B \right) \cdot \zs\right),f \right]_{\#} 
									&\textrm{ on } \mathbb{R}_x^3\times\mathbb{R}^3_p\times[0 ,T] \\[2pt]
									\disp i\hbar \dpp{h}{t} =  \left[ \left( \frac{p^2}{2m} +U  +\hbar\left( p \wedge K-   B \right) \cdot \zs  \right),h \right]_{\#} &\textrm{ on } \mathbb{R}_x^3\times\mathbb{R}^3_p\times[0 ,T]\\[2pt]
									\left.f \right|_{t=0} =\mathcal{W}^\hbar [\zy_0]  & \textrm{ in } \mathbb{R}_x^3\times\mathbb{R}_p^3 \\[2pt]
									\left.h \right|_{t=T}  =f_T  & \textrm{ in }\mathbb{R}_x^3\times\mathbb{R}_p^3 \\[2pt]
									\disp 	\zg'  \dst{u}{t}-	\zg u   =-   \frac{1}{2} \textrm{tr} \int_{\mathbb{R}_x^3\times\mathbb{R}_p^3}   \Theta^-_{\nabla_u U} [h] \;f    \dif x  \dif p&\textrm{ on } [0,T]
								\end{array}\right.\label{Opt_cond_Q}
							\end{align}
							
							\subsection{OC problem for the classical deterministic dynamics}
							
							This work aims to clarify the mathematical formulation of the classical limit associated with the OC problem \eqref{Opt_cond_Q}. We expect that as $\hbar$ goes to zero, the solution tends to localize, and the quantum dynamics degenerates to the deterministic classical transport described by a single trajectory in the phase-space. The solution may be interpreted as the evolution of a spin vector traveling in the phase-space and, at the same time, rotating around to the total magnetic field given by the sum of the external and the Rashba effective fields. In this section, we provide a formal derivation of the classical limit of the OC problem, and we discuss some of the difficulties emerging in the comparison between the quantum and deterministic OC problems. 
							We start by considering the formal classical limit of the quantum dynamics resulting from retaining the leading terms of the expansion given in Eq. \eqref{form_hbar_exp_fk}. The $\hbar=0$ terms provide the Liouville-Rashba equation for the classical particle density $f: \mathbb{R}_x^3\times \mathbb{R}_p^3\times [0,T]\rightarrow\mathbb{R}^3 $
							\begin{align}
								&\dpp{f}{t}+\frac{p}{m} \cdot \nabla_{x} f +E  \cdot \nabla_{p} f + \left( B - p\wedge K \right) \wedge f =0&\textrm{ on } \mathbb{R}_x^3\times\mathbb{R}_p^3\times[0 ,T] \;.\label{Liou_f_1D}
							\end{align}
							We have defined the electric field $E\doteq -  \nabla_{x}U$.
							The characteristic curves associated to the classical Liouville-Rasha equation can be interpreted as the phase-space trajectories followed by non interacting particles. Each particle is described by the phase-space coordinates $(x,p)\in \mathbb{R}^3_x\times \mathbb{R}^3_p$, and one vector $d\in\mathbb{R}^3 $, representing the direction of the spin or dipole momentum. The single-particle dynamics is expressed by the Hamiltonian system 
							\begin{align}
								\left\{
								\begin{array}{ll}
									\disp \dot{x} = \frac{p}{m} & \textrm{on } t\in [0,T]\\[2mm]
									\disp\dot{p} =	E(x,u)& \textrm{on } t\in [0,T]\\[2mm]
									\dot{d}=2    ( p\wedge K - B )\wedge  d& \textrm{on } t\in [0,T]\\[2mm]
									x= \overline{x} ,\; p=\overline{p},\; d=\overline{d} & \textrm{in } t=0
								\end{array}
								\right. . \label{Char_traj_1D}
							\end{align}
							Here, $\overline{x}$, $\overline{p}$ and $\overline{d}$ denote the initial position, momentum and spin, respectively. 
							We now set the control problem at the single-particle level. Consistently with the OC strategy adopted in the quantum dynamic case, we derive the OC problem for the classical evolution and we adopt as a goal the necessity to steer the particle towards the phase-space coordinates $(x_T,p_T)$ and align the spin direction to the desired vector $d_T$, at the final time $T$. In close analogy to Eq. \eqref{goal_funct}, we define the following goal functional, that we would like to minimize
							\begin{align}
								\Phi'(x,p,d)=& \frac{\nu_x}{2}|x(T) -x_T|^2+\frac{\nu_p}{2} p^2(T)- \nu_d\; d(T)\cdot d_T.\label{goal_funct_sing_tj}
							\end{align} 
							The cost of the control is given by Eq. \eqref{cost_func} as before. 
							We are led to consider the OC problem for the deterministic classical dynamics
							\begin{align}
								\text{Cl-OC: }\: \underset{ u }{\text{min}} \; \left[\Phi'(x,p,d)+k(u)\right] &\quad  \textrm{s. t. } \quad \textrm{Eq. } \eqref{Char_traj_1D} \textrm{ holds true.}\label{OP_traj_1D}
							\end{align}
							The optimality system is obtained by varying the classical Lagrangian with respect to the set of variables $(x,p,d,x^h,p^h,d^h,u)$ 
							\begin{align*}
								\mathcal{L}=& \int_0^T \left[\left( 	\dot{x} - \frac{p}{m} \right) p^{h} +  	\left(\dot{p} -E \right)x^{h}+ \left(	\dot{d} -2    ( p \wedge K - B )\wedge  d \right)\cdot p_d^{h}\right]\dif t+ k + \Phi'\;. 
							\end{align*}
							In this single particle formulation, the role of the adjoint function $h$ is played by the adjoint trajectory $p_d^h$.  
							Standard calculations provide the following set of adjoint equations 
							\begin{align}
								\left\{
								\begin{array}{ll}
									\disp \dot{x}^h=- \frac{p^h}{m}   +2  K\wedge \left(	 p_d^{h}   \wedge    d  \right) & \textrm{on } t\in [0,T]\\[2mm]
									\disp\dot{p}^h  =-\dpp{E}{x} x^h +   2  \left( p \wedge\dpp{K}{x} - \dpp{B}{x}   \right)  \cdot \left(	 p_d^{h}   \wedge    d  \right)& \textrm{on } t\in [0,T]\\[2mm]
									\disp\dot{p}_d^h =2    ( p\wedge K - B ) \wedge p_d^h& \textrm{on } t\in [0,T]\\[2mm]
									\disp x^h(T)=  - \nu_p p(T),  \;  
									\disp p^h(T)= \nu_x(x_T-x(T)),	 \;\disp p_d^h(T)=   \nu_{d} d_T   &\textrm{in } t=T
								\end{array}
								\right.\;.
							\end{align}
							It is convenient to define $\zh^h =  p_d^{h}   \wedge    d$ and we obtain
							\begin{align}
								\left\{
								\begin{array}{ll}
									\disp \dot{x}^h=- \frac{p^h}{m}   +2  K\wedge \zh & \textrm{on } t\in [0,T]\\[2mm]
									\disp\dot{p}^h  =-\dpp{E}{x} x^h +   2  \left( p \wedge\dpp{K}{x} - \dpp{B}{x}   \right)  \cdot  \zh^h & \textrm{on } t\in [0,T]\\[2mm]
									\disp\dot{\zh^h} =2    ( p\wedge K - B ) \wedge \zh^h& \textrm{on } t\in [0,T]\\[2mm]
									\disp x^h(T)=  - \nu_p p(T)   ,\;
									p^h(T)= \nu_x(x_T-x(T)),\;
									\zh^h (T)=   \nu_{d} d_T \wedge d(T)  &\textrm{in } t=T
								\end{array}
								\right.\;. \label{char_adj_1D_03}
							\end{align}
							The equation for the control is given by
							\begin{align}
								\label{Opt_cont_conc_lim_u}
								&\zg'   \dst{u_i}{t}-	\gamma u_i=-x^h\frac{\partial E}{\partial u_i}& i=1,\ldots,d.
							\end{align}
							We have obtained the following optimality system associated to the OC for the deterministic classical transport
							\begin{align}
								\textrm{Cl-OS: }\;\left\{
								\begin{array}{ll}
									(x,p,d) &   \textrm{Classical dynamics, Eq. } \eqref{Char_traj_1D} \\
									(x^h,p^h,\zh^h)&  \textrm{Adjoint variables, }\textrm{Eq. } \eqref{char_adj_1D_03} \\
									u &\textrm{Control,  }  \textrm{Eq. } \eqref{Opt_cont_conc_lim_u}  
								\end{array}
								\right.\label{cl_OC_prob}
							\end{align}
							We show that in the limit of small $\hbar$ the quantum OC problem converges to the deterministic single particle OC problem. Comparing the two sets of optimality conditions, we see that the main difference concerns the equation for the control which is expressed in the two systems by physical and adjoint variables, which have different mathematical formulation and whose relationship is not straightforward. 
							
							\section{Classical deterministic limit of the quantum OC problem}\label{Sec_cl_lim}

							We show that in the limit $\hbar\rightarrow 0$, the OC problem \eqref{Opt_cond_Q} for quantum dynamics converges to the analogous deterministic OC problem \eqref{char_adj_1D_03} associated with the classical evolution of the particle. 
							In order to make such a comparison possible, it is necessary to specify in which way the particle wave function scales with respect to the Planck constant. When $\hbar$ goes to zero, we expect that the be Broglie length of the quantum particle vanishes and the particle mass concentrates in a single point. This limit can be investigated by considering coherent states. For any L2 wave function $\zf$, we define the coherent state associated with the phase-space coordinates $(\overline{x},\overline{p})$, the wave function $\zy^\hbar= \frac{1}{\hbar^{3/4}}\zf\left(\frac{x-\overline{x}}{\sqrt{\hbar}}\right)e^{i\frac{\overline{p}}{\hbar}x}$. Hereafter, for simplicity, we will assume that the initial wave function of the system is given by the tensor product of a scalar wave function and a spinor $\chi\in \mathbb{C}^2$, with $|\chi|^2=\chi_1^2+\chi_2^2=1$, i. e. $\zf=\zf_o\otimes\chi $, with $\zf_o\in L^2\left(\mathbb{R}^3_x\right)$. With this choice, the spin expectation value of the coherent state is independent of $\hbar$ and is given by $d_i=\int_{\mathbb{R}^3_x} \langle\zy^\hbar,\zs_i\zy^\hbar \rangle_{\mathbb{C}^2}\dif x = \langle\chi,\zs_i\chi\rangle_{\mathbb{C}^2} $, with $i=1,2,3$. 
							The Wigner transform of the coherent state $f^\hbar(t) = \mathcal{W}^\hbar[\zy^\hbar(t)]$ converges as $\hbar$ goes to zero, to the Dirac measure $
							\zd\left(x-x(t)\right)\zd\left(p-p_{c}(t)\right) d(t)$,
							where the classical phase-space trajectories and the spin vector are obtained as a solution of Eq. \eqref{Char_traj_1D}. With these considerations, we can state our main result
							\begin{theorem}\label{Theor_lim_OC}
								Under AS, for any set of nonnegative control widths $\nu_x,\nu_p,\zg,\zg'$, goals $(x_T,p_T,d_T)$, phase-space coordinates $(\overline{x},\overline{p})$, and quantum states $\zf\in L^2(\mathbb{R}^3,\mathbb{C}^2)$, with $\zf=\zf_o\otimes\chi $ and $\chi\in \mathbb{C}^2$, the solution $(f^\hbar,u^\hbar)$ of the OC problem \eqref{Opt_cond_Q} with initial condition $f^\hbar(0) =  \mathcal{W}^{\hbar}[\zf_{(\overline{x},\overline{p})}^{\hbar}]\in H^1_p$, 
								after possible extraction of subsequence,  converges $(f^\hbar,u^\hbar) \stackrel{\hbar\rightarrow 0}{\longrightarrow} (f^0,u^0)$ as distributions, where $f^0=\zd\left(x-x(t)\right)\zd\left(p-p(t)\right) d(t)$ and the functions $(u^0,x,p,d)$ are a solution of the deterministic classical OC problem of Eq. \eqref{cl_OC_prob}, with initial conditions $(\overline{x},\overline{p},\overline{d})$, where $\overline{d}=\langle\chi,\zs\chi\rangle_{\mathbb{C}^2}$.
							\end{theorem}
							
							We split the proof of Th. \ref{Theor_lim_OC} in  three main steps: $I:$ Study of the classical limit of the Wigner dynamics; $II:$ Study of the analogous problem for the adjoint function; $III:$ Limit $\hbar\rightarrow 0$ of the optimality system and comparison with the deterministic dynamics. We briefly discuss the various points. $I:$ First, we investigate the concentration limit of the Wigner pseudo-distribution function when $\hbar$ goes to zero. The well-posedness of the dynamics for $\hbar =0$ has already been investigated in Sec. \ref{Sec_model}, Th. \ref{teor_wellpos_W_dyn}. We focus on the initial datum of the Wigner matrix function. In the limit $\hbar \rightarrow 0$ the initial condition of the Wigner dynamics is given in the space of distributions. The initial condition of the dynamical problem becomes a delta-like distribution in phase-space. The classical limit of the initial datum is investigated by considering coherent states associated with phase-space translations and scaling of L2 wave functions. This is done in Lemma \ref{Lemma_lim_meas_IC}. The proof relies uniquely on the definition of the Wigner transformation.  In Th. \ref{Theor_lim_eq_Wig_f}, we show that the resulting Wigner measure propagates along the trajectories associated with the classical Hamitlonian of a particle with spin in the presence of the Rashba-Zeeman field. $II:$ After such a preliminary step, we perform a similar analysis for the adjoint function. The main difference from the previous case concerns the formulation of the boundary value of the equations and the regularity of the solution. The adjoint function $h$ cannot be seen as the Wigner transform associated with any Schr\"odigner wave function. The adjoint function reveals to be more regular than the Wigner function, so that the classical limit is taken in the space of continuous functions instead of in the distribution space. We proceed as follows: We introduce a cut-off in the phase-space in order to ensure the L2 and H1 integrability of the solution. We study the classical limit and analyze the regularity of the limit solution. This is done in Th. \ref{teor_wellpos_h_dyn}. In the final step $III$, we show that the control can be determined by integrating along a single trajectory of a Hamiltonian system. Finally, we remove the cut-off introduced previously.  %
							\subsubsection*{$I:$ Study of the classical limit of the Wigner dynamics}
							The quantum state of the system is represented by a wavefunction $\zy^\hbar$ which is a solution of the Schr\"odinger problem
							\begin{align*}
								\left\{
								\begin{array}{ll}
									i\hbar \dpp{\zy^\hbar}{t} =\textrm{Op} \left(\mathcal{H}\right) \zy^\hbar &\textrm{on } \mathbb{R}^3_x\\
									\zy(0)=\zf_{(\overline{x},\overline{p})}^\hbar& \textrm{in } t=0
								\end{array}
								\right.\;.
							\end{align*}
							The concentration of the solution toward a single point in phase-space emerges after suitable scaling of the initial condition of the wave function with respect to $\hbar$. We consider the following family of coherent states associated with any function $\zf\in L^2 \mathbb{R}^3,\mathbb{C}^2)$
							\begin{align}
								\zf_{(\overline{x},\overline{p})}^\hbar\doteq \frac{1}{\hbar^{3/4}}\zf\left(\frac{x-\overline{x}}{\sqrt{\hbar}}\right)e^{i\frac{\overline{p}}{\hbar}x}\,,\label{scaled_coh_wf}
							\end{align}
							obtained by translating $\zf$ into the phase-space by the coordinates $(\overline{x},\overline{p})$. The behavior of the Wigner transform of $\zy^\hbar$ is discussed in Lemma \ref{Lemma_lim_meas_IC}.
							As indicated in the previous section, the Wigner matrix-function  $f^\hbar=\mathcal{W}^\hbar[\zy]$ is a solution of the Wigner equation
							\begin{align*}
								\left\{
								\begin{array}{ll}
									i\hbar \dpp{f^\hbar}{t} = \left[\mathcal{H} , f^\hbar \right]_\#&\\
									f^\hbar(0)=\mathcal{W}^\hbar[\zy_0]& \mathrm{in } \; t=0
								\end{array}
								\right.\;. 
							\end{align*}
							The conservation of the L2 norm of the wave function, has as a consequence that the sequence $f^\hbar$ forms a set of distributions uniformly bounded. The idea of the proof is based on the following straightforward estimate
							\begin{align*}
								\left|\left\langle f^\hbar ,h \right\rangle_{\mathrm{HS}} \right| \leq C \| \zy^\hbar \|_{L^2(\mathbb{R}_x^3)}  \| \mathcal{F}_{v\rightarrow\zh} h\|_{C^0(\mathbb{R}_x^3, L^1 (\mathbb{R}_\zh^3))}\;,
							\end{align*}
							where $h$ is a test function. 
							As a well-known consequence, after the possible extraction of subsequences, the $f^\hbar$ distributions have $*-w$ limit. A more detailed analysis of the Wigner transform associated with the scaled coherent states is given in the following 
							\begin{lemma}\label{Lemma_lim_meas_IC}
								For any $\overline{x},\overline{p} \in \mathbb{R}^3$, define $\zf_{(\overline{x},\overline{p})}^\hbar $ according to Eq. \eqref{scaled_coh_wf}, where $\zf=\zf_o\otimes\chi $, $\chi\in \mathbb{C}^2$, and $\|\zf\|_{L^2(\mathbb{R}^3_x,\mathbb{C}^2)}=1$. Let $\{\ze_n\}_{n\in\mathbb{N}} $ be any monotonic decreasing sequence that converges to zero, and $\zf^{\ze_n}(t) = e^{-\frac{i}{\ze_i}\mathrm{Op}(\mathcal{H}) t}\zf^{\ze_n}_{(\overline{x},\overline{p})} $ be the solution of the Schr\"odinger equation with initial condition $\zf^{\ze_n}_{(\overline{x},\overline{p})}$. After possible extraction of an appropriate subsequence, the Wigner transform $f^n \doteq  \mathcal{W}^{\ze_n}[\zf_{(\overline{x},\overline{p})}^{\ze_n}]$ is weak-$*$ convergent in the space of distribution $\mathcal{S}'\left(\mathbb{R}_x^3\times \mathbb{R}_p^3,\mathbb{C}^2\times \mathbb{C}^2\right)$. 
								In particular, for $t=0$ we have  
								\begin{align}
									&\lim_{\hbar \rightarrow 0} \textrm{tr} \left\{f^\hbar (0) \zs \right\} =  \overline{d}\; \zd\left(x-\overline{x}\right)\zd\left(p-\overline{p}\right) \label{lim_fh}\\
									&\lim_{\hbar \rightarrow 0} \textrm{tr} \left\{f^\hbar (0) \right\} = \zd\left(x-\overline{x}\right)\zd\left(p-\overline{p}\right)\;,\label{lim_fh0}
								\end{align}
								where $ \overline{d}=\langle\chi,\zs\chi\rangle_{\mathbb{C}^2}$  is the expectation value of the spin.
							\end{lemma}
							\proof The proof is standard, except for the presence of the spin degree of freedom. Assuming $\mathsf{h}\in \mathcal{S}(\mathbb{R}^3_x\times\mathbb{R}^3_p,\mathbb{C}^4)$, straightforward calculations show that
							\begin{align*}
								\left\langle f^\hbar ,\mathsf{h} \right\rangle =&\frac{1}{(2\pi)^{3/2}}	\sum_{i,j}\int_{\mathbb{R}^{3}_x\times\mathbb{R}^{3}_\zh} \zf_{i} \left(x+\frac{ \sqrt{\hbar}\zh}{2} \right) \overline{\zf_{j} \left(x-\frac{\sqrt{\hbar}{ \zh}}{2}\right)} (\mathcal{F}_{p\rightarrow \zh}\mathsf{h}_{ji} )\left(x\sqrt{\hbar}+\overline{x},\zh \right)  e^{i\overline{p} \zh}\dif { \zh}\dif x \;.
							\end{align*}
							We note that 
							\begin{align*}
								\int_{\mathbb{R}^{3}_x\times\mathbb{R}^{3}_\zh} \left|\zf_{i} \left(x+\frac{ \sqrt{\hbar}\zh}{2} \right) \right| \left|{\zf_{j} \left(x-\frac{\sqrt{\hbar}{ \zh}}{2}\right)}\right| \left|(\mathcal{F}_{p\rightarrow \zh}\mathsf{h}_{ji} )\left(x\sqrt{\hbar}+\overline{x},\zh \right) \right| \dif { \zh}\dif x \\
								\leq 
								\|\zf_i\|_{L^2(\mathbb{R}^3_x)}\|\zf_j\|_{L^2(\mathbb{R}^3_x)}\int_{ \mathbb{R}^{3}_\zh} \sup_{x\in \mathbb{R}^3_x} \left|(\mathcal{F}_{p\rightarrow \zh}\mathsf{h}_{ji} )\left(x,\zh \right) \right| \dif { \zh} < \infty\;.
							\end{align*}
							Lebesgue dominated convergence theorem gives
							\begin{align*}
								\lim_{\hbar\rightarrow 0 }\left\langle f^\hbar ,\mathsf{h} \right\rangle = \sum_{i,j}\mathsf{h}_{ji} \left(\overline{x},\overline{p} \right) \int_{\mathbb{R}^{3}_x\times\mathbb{R}^{3}_\zh} \zf_{i} \left(x \right) \overline{\zf_{j} \left(x\right)}  \dif x \;.
							\end{align*}
							Taking $\mathsf{h}=\zs \mathsf{h}'$ or $\mathsf{h}=\zs_0 \mathsf{h}'$ with $\mathsf{h}'\in \mathcal{S}(\mathbb{R}^3_x\times\mathbb{R}^3_p)$, we obtain Eq. \eqref{lim_fh} and Eq. \eqref{lim_fh0}, respectively.
							Concerning the first part of the Lemma, 
							from Lemma \eqref{Lemma_ess_sa_H} we have that the Hamiltonian $\mathcal{H}$ generates a strongly continuous semi-group. The $\{\zf^{\ze_n}(t)\}_{n\in \mathbb{N}}$ is a well-defined sequence in $L^2(\mathbb{R}^d_x,\mathbb{C}^2)$ with a unit norm.
							Applying standard arguments \cite{Lions_93}, we verify that the sequence of the associated Wigner distributions $\{f^{n}= \mathcal{W}^{\ze_n}[\zf^{\ze_n}(t)]\}_{n\in \mathbb{N}}$   is uniformly bounded and thus compact in the weak$-*$ topology. 
							\endproof
							To make the analysis of the limit of the Wigner function more consistent, it would be desirable to avoid any reference to the original Schr\"odinger wavefunction description, and derive the previous results using, for example, the conservation properties of the L2 norm of the solution of the Wigner equation, or alternative estimates which can be proven inside the Wigner formalism. Since, in the limit, $f$ concentrates towards a delta, alternative approaches that rely only on the Wigner formalism seem challenging. On the other hand, different choices on the initial conditions may lead to more regular limits. For initial data uniformly into L2, estimations on the norm growth rates are sufficient to guarantee the existence of the limit toward the solution of a Liouville equation without any concentration effect. This will be the strategy adopted for the adjoint equation in Sec. \ref{sec_adj_eq} below. 
							
							
							The strong formulation of Wigner dynamics is expressed in compact form using the Moyal formalism by Eq. \eqref{Wig_dyn_Moy}, $i\hbar \dpp{f}{t} = \left[\mathcal{H}_0 + \mathcal{H}^K,f \right]_{\#} $, or explicitly, by Eq. \eqref{Wig_f0_fk_01}. Since in the classical limit the quantum spin dynamics is expressed in a distributional sense, we pass to the weak formulation of Eq. \eqref{Wig_dyn_Moy}. The scalar product of the natural Hilbert space associated with the Wigner representation of the quantum dynamics of particles with spin, is given by $\langle f,h \rangle \doteq\frac{1}{2} \mathrm{tr}\int_{\mathbb{R}^{6}_{x,p} } f^\dag h \dif x \dif p$.  We consider the expression $i\hbar \left\langle \dpp{f}{t}, h \right\rangle  =\left\langle \left[\mathcal{H}_0 + \mathcal{H}^K,f \right]_{\#} , h \right\rangle $ where $h$ is any function belonging to $\mathrm{HS}$. In Lemma \ref{lemma_anti_symm_H} we show that the Moyal product is anti-symmetric i. e.
							$		\left\langle 	[\mathcal{H} ,f]_\#  ,h \right\rangle=-\left\langle f,	[\mathcal{H}, h ]_\#   \right\rangle$, and the Wigner dynamics is expressed in weak form by $i\hbar \left\langle f,\dpp{h}{t}  \right\rangle  =\left\langle f,\left[\mathcal{H},h \right]_{\#} \right\rangle $, with $h\in \mathrm{HS}$. The explicit form of this expression is given by Eq. \eqref{Wig_f0_w_01}.
							In order to avoid confusion, we adopt the following notation. The symbol $f^\hbar = \mathcal{W}^\hbar[\zf]$ denotes the Wigner transform with values into the space of $2\times2$ hermitian matrices. The Pauli components along the coordinate axes are indicated by an arrow ${\vec{f}}_i=\frac{1}{2} \mathrm{tr}\left\{ f \zs_i \right\} $ with $i=1,2,3$, and the component along the identity matrix by the subscript zero ${f}_0=\frac{1}{2} \mathrm{tr}\left\{ f  \right\}$. 
							The following theorem shows that in the classical limit the dynamics can be approximated by a system of classical Liouville equations coupled via the Rashba-Zeeman fields and have as initial datum point sources. 
							\begin{theorem}\label{Theor_lim_eq_Wig_f}
								Under AS, in the limit $\hbar \rightarrow 0$, any solution $f^\hbar$ of Eq. \eqref{Wig_f0_fk_01} with initial condition $f^\hbar (0)= \mathcal{W}^\hbar[\zf_{(\overline{x},\overline{p})}^\hbar]$ and $\zf=\zf_o\otimes \chi$, converges to a weak solution $f$ of the associated classical Liouville-Rashba-Zeeman equation
								\begin{align}
									\left\{\begin{array}{lll}
										\disp	\dpp{f_0}{t}  +  \frac{p}{m}\cdot\nabla_xf_0 -\nabla_xU \cdot\nabla_p f_0 =0 &\textrm{on } \mathbb{R}_x^3\times\mathbb{R}^3_p\times[0 ,T] \\[6pt]
										\disp	\dpp{\vec{f}}{t} + \frac{p}{m}\cdot\nabla_x \vec{f} -\nabla_xU \cdot\nabla_p \vec{f} = 2    ( p\wedge K - B  )\wedge \vec{f} &\textrm{ on } \mathbb{R}_x^3\times\mathbb{R}_p^3\times[0 ,T] \\[6pt]
										\disp  \left. \vec{f}\; \right|_{t=0} = \overline{d}\; \zd\left(x-\overline{x}\right)\zd\left(p-\overline{p}\right) &\textrm{ on } \mathbb{R}_x^3\times\mathbb{R}_p^3 \\
										\left. f_0\right|_{t=0}=  \zd\left(x-\overline{x}\right)\zd\left(p-\overline{p}\right) &\textrm{ on } \mathbb{R}_x^3\times\mathbb{R}_p^3 
									\end{array}
									\right.\label{Liouv_Rash_01}
								\end{align}
								where   $\overline{d} =\langle\chi,\zs\chi\rangle_{\mathbb{C}^2} $. The solution of the second equation is given by $
								\vec{f}=  d(t) \zd\left(x-x(t)\right)\zd\left(p-p(t)\right)$, 
								where the phase-space coordinates  $(x(t),p(t))$ and the spin $d$ are obtained integrating the Hamitonian system of Eq. \eqref{Char_traj_1D}.
							\end{theorem}
							\proof
							As a distribution, a weak solution of Eq. \eqref{Wig_f0_fk_01} $, f\in \mathcal{S}'(\mathbb{R}^3_x\times\mathbb{R}^3_p,\mathbb{C}^4)$ satisfies  
									\begin{align}
										\int_0^\infty  \left\langle f, \left(\dpp{\mathsf{h}}{t}  + \frac{p}{m} \cdot\nabla_x \mathsf{h}   - \Theta^-_U [\mathsf{h}] \right)  \right\rangle \dif t =&-	\sum_{k=1}^3  \int_0^\infty  \left\langle f,   \zs_k  \left( \sum_{i,j=1}^3   \ze_{ijk} \Theta^+_{K_{j}}[p_{i} \mathsf{h} ]  +\Theta^+_{ B_k}[\mathsf{h}]   \right)   \right\rangle \dif  t  \nn  \\
										&+ \int_0^\infty \left(   \mathrm{Re }  \left\langle f , Q [\mathsf{h}]   \right\rangle   -\mathrm{Im } \left\langle f , D [\mathsf{h}]   \right\rangle    \right)   \dif t  \; , \nn\\
										&\forall\mathsf{h}\in \mathcal{S}(\mathbb{R}^3_x\times\mathbb{R}^3_p,\mathbb{C}^4), \label{Wig_f0_w_01}
									\end{align}
									where we have defined
									\begin{align*}
										D[\mathsf{h}]	\doteq &-\sum_{i,j,k=1}^3  \zs_k  \ze_{ijk}\frac{\hbar^2}{2}   \dpp{}{x_{i}}\Theta^-_{K_{j}}\left[\mathsf{h} \right] \\
										Q [\mathsf{h}]	\doteq & \hbar \sum_{k=1}^3 \zs_k \left[\sum_{i,j=1}^3 \ze_{ijk}    \left( 	\Theta^-_{K_j}[ p_i \mathsf{h}]+ \frac{1}{2}   	\Theta^+_{K_j}\left[\dpp{\mathsf{h} }{x_i}\right]\right)  -  \Theta^-_{ B_k}[ \mathsf{h}]\right] \; .
									\end{align*}
									The details of the derivation of Eq. \eqref{Wig_f0_w_01} are given in Appendix \ref{app_der_w_form}. 
									In the limit $\hbar\rightarrow 0$ the last two terms in the right side of Eq. \eqref{Wig_f0_w_01} vanish and the remaining pseudo-differential operators simplify to first-order derivative or multiplication operators. Assuming $\mathsf{h}\in \mathcal{S}(\mathbb{R}^3_x\times\mathbb{R}^3_p,\mathbb{C}^4)$, 
									and using $		\zd_+ K_j\left(x , \zh \right) = 2 K (x) + \frac{\hbar}{2}\int_{-1}^1 \textrm{sgn}(s)\; \zh\cdot \nabla_x 	K_j \left(x +\frac{\hbar \zh}{2}  s \right)\dif s
									$ (see Eq. \eqref{delta_k_+} for the definition of $\zd_+$), we show that $\zs_k  \Theta^+_{K_{j}}[p_{i} \cdot  ] \stackrel{*}{\rightharpoonup} 2\zs_k K_{j}p_{i}  $. 
									We have
									\begin{align*}
										\left\langle  f^\hbar ,     \zs_k  \left(\Theta^+_{K_{j}}[p_{i} \mathsf{h} ]- 2K_{j}p_{i}\mathsf{h}   \right)
										\right\rangle  =& \frac{ 1}{(2\pi)^{3}} \mathrm{tr} \;  \int_{\mathbb{R}^{12}} f^\hbar(x,p) \left[\zd_+ K_{j}\left(x , \zh \right)p'_{i} - 2K_{j}p_{i} \right] \zs_k  \mathsf{h}  \left(x,p' \right) e^{ i\zh(p -p') } \mathrm{d}  \zh  \mathrm{d}  p'\mathrm{d}x\mathrm{d}   p\\
										=& \frac{ i\hbar }{ (2\pi)^{3/2} } \mathrm{tr} \;  \int_{\mathbb{R}^{6}}\zr^\hbar  \left(x+\frac{\hbar{ \zh}}{2},x-\frac{\hbar{ \zh}}{2} \right) \left[\zh\cdot H^\hbar(x,\zh) \right]\zs_k \\
										& \hspace{150pt}\times\dpp{}{\zh_{i'}}\mathcal{F}_{p'\rightarrow\zh}\left( \mathsf{h}  \left(x,p' \right)  \right) \dif \zh   \dif x\;,
									\end{align*}
									where  we have defined $H^\hbar(x,\zh)\doteq \int_{-1}^1 \textrm{sgn}(s)\;  \nabla_x 	K_j \left(x +\frac{\hbar \zh}{2}  s \right)\dif s$.
									We obtain the estimate
									\begin{align*}
										\left|\left\langle  f^\hbar ,  \zs_k  \left(\Theta^+_{K_{j}}[p_{i} \mathsf{h} ]- 2K_{j}p_{i} \mathsf{h}  \right)
										\right\rangle \right| 
										\leq&  \frac{ \hbar \|K\|_{ W^{1,\infty}(\mathbb{R}^3_x,\mathbb{R}^3)}}{  (2\pi)^{3/2} }   \|\zf\|_{L^2(\mathbb{R}_x^3,\mathbb{C}^2)}^2  \max_{r,s}  \int_{\mathbb{R}^{3}}\sup_{x}|\zh| \left|\dpp{\mathcal{F}_{p\rightarrow\zh}\mathsf{h}_{rs}}{\zh_{i}} \left(x,\zh \right) \right|    \dif \zh.
									\end{align*}
									Similarly, standard estimates provide the limits $ \Theta^-_{U}[ \cdot ]
									\stackrel{*}{\rightharpoonup}  	 \nabla_{x}U\cdot  \nabla_p $ and $
									\zs_k \Theta^+_{ B_k}[\cdot]    \stackrel{*}{\rightharpoonup}\zs_k 2  B_k $. The remaining terms go to zero for small $\hbar$. As an example, we consider
									\begin{align*}
										\left\langle  f^\hbar ,   \zs_k    \dpp{}{x_{i}}\Theta^-_{K_{j}}\left[\mathsf{h} \right]
										\right\rangle  
										=&    \frac{ i}{ \hbar(2\pi)^{3/2} } \mathrm{tr} \;  \int_{\mathbb{R}^{6}}\zr^\hbar  \left(x+\frac{\hbar{ \zh}}{2},x-\frac{\hbar{ \zh}}{2} \right) \left[\dpp{}{x_{i}}\zd_- K_{j}\left(x , \zh \right) \right]\zs_k (\mathcal{F}_{p\rightarrow\zh}   h)  \left(x,\zh  \right) \dif \zh   \dif x\;.
									\end{align*}
									Using $\dpp{}{x_{i}}	\zd_- K_j\left(x , \zh \right)= \frac{\hbar}{2}\int_{-1}^1 \zh\cdot \nabla_x \dpp{}{x_{i}}	K_j \left(x +s\frac{\hbar \zh}{2}  \right)\dif s$, we obtain the estimate
									\begin{align*}
										\left|\left\langle  f^\hbar ,   \zs_k    \dpp{}{x_{i}}\Theta^-_{K_{j}}\left[\mathsf{h} \right]
										\right\rangle  \right| 
										\leq& \frac{ \hbar^2 \|K\|_{ W^{2,\infty}(\mathbb{R}^3_x,\mathbb{R}^3)}}{  2(2\pi)^{3/2} }   \|\zf\|_{L^2(\mathbb{R}_x^3,\mathbb{C}^2)}^2  \max_{r,s}  \int_{\mathbb{R}^{3}}\sup_{x}|\zh| \left| \mathcal{F}_{p\rightarrow\zh} \mathsf{h}_{rs} \left(x,\zh \right) \right|    \dif \zh\;.
									\end{align*}
									The initial conditions follow from Eq. \eqref{lim_fh} and Eq. \eqref{lim_fh0} in Lemma \ref{Lemma_lim_meas_IC}. Under AS the characteristic curves of the Liouville-Rashba system exist globally in time \cite{Reed_Simon_vol_II},
									and the characteristics of Eq. \eqref{Liouv_Rash_01} are given by the system \eqref{Char_traj_1D}. 
									%
									\hfill $\square$\\
									\subsection*{$II$ Adjoint problem}\label{sec_adj_eq}
									We pass to the study of the adjoint problem. The analysis proceeds similarly to the case of the Wigner function, with the important difference that the adjoint function $h$ is no longer the Wigner transform of some wavefunction. This difference appears clearly in the choice of boundary values. For the Wigner problem, the initial value is fixed by physical considerations and coincides with the Wigner transform of the particle wavefunction. Concerning the adjoint equation, the boundary value is fixed at the final time $T$ and represents the distance to the target function associated with the minimization problem. In order to have compatibility with the classical deterministic case, the distance function is chosen to be quadratic in the position and in the momentum coordinates. In the limit $\hbar\rightarrow 0 $, this choice leads to the usual Euclidean norm of the position-momentum error. As a consequence, the final value of the adjoint equation does not belong to the L2 space. 
									In order to ensure bonds on the L2 norm of the solution, we introduce a cut-off on the final condition. Our strategy is based on the following considerations. 
									Unlike the Wigner function, in the limit $\hbar \rightarrow 0 $ the adjoint function does not concentrate in the phase-space. In fact, the adjoint function represents the evolution of the distance to the target, which is a smooth function. In order to study the limit of the OC quantum problem toward the classical deterministic OC, it is sufficient to ensure that the application of the Wigner measure on the adjoint function is well defined. The reason is that the optimality equation for the control parameters depends on the integral of the product between the Wigner and the adjoint function (see the last equation of the OC problem \eqref{Opt_cond_Q}). Natural considerations based on the study of the Hamiltonian trajectories ensure that the limit single-particle trajectories are bounded in the phase-space. As a consequence, in the limit, the adjoint function will be tested within a bounded region of the phase-space. 
									For any initial value of the physical position of the deterministic particle to be controlled, it is always possible to fix a cut-off sufficiently big to contain the support of the Wigner measure. 
									
									The conservation of the L2 norm of the adjoint function follows from the anti-symmetry of the Hamiltonian operator. We have 
									\begin{lemma}\label{lemma_anti_symm_H}
										The Moyal commutator of the Hamiltonian symbol is anti-symmetric with respect to the duality $\langle f,h \rangle=\frac{1}{2}\mathrm{tr}\int_{\mathbb{R}^{6}_{x,p} } f^\dag h \dif x \dif p $. For any $f,h \in \mathrm{HS}$ 
										we have
										\begin{align}
											\left\langle 	[\mathcal{H} ,f]_\#  ,h \right\rangle=-\left\langle f,	[\mathcal{H} ,h ]_\#   \right\rangle\;.\label{symm_H}
										\end{align}
									\end{lemma}
									\proof
									We limit ourselves to discussing the result for the Rashba Hamiltonian, the anti-symmetry of the other terms can be easily verified. First, we note that for any $i,j=1,2,3$, we have
									\begin{align}
										&\int_{\mathbb{R}^{6}}  \left[ p_iK_j, f\right]_\# h \dif p\dif x  	  =- \int_{\mathbb{R}^{6}}  f \left[ p_iK_j, h\right]_\#  \dif p\dif x  	\label{Wig_sym_01} \\
										&\int_{\mathbb{R}^{6}}  \left\{p_iK, f\right\}_\# h \dif p\dif x  	  = \int_{\mathbb{R}^{6}}  f \left\{ p_iK, h\right\}_\#  \dif p\dif x  	\;.\label{Wig_sym_02}
									\end{align}
									We verify Eq. \eqref{Wig_sym_01}. From the definition of the Moyal product,  we have
									\begin{align}
										\left[  p_i K_j , f \right]_{\#} =&
										i\hbar p_i  \Theta^-_{K_j}[f]   - i \frac{\hbar}{2}  \Theta^+_{K_j}\left[\dpp{f}{x_i}\right]  \;.	\label{moy_prod_pKf_01}
									\end{align}
									Using integration by part, the first term of Eq. 	\eqref{moy_prod_pKf_01} can be written as 
									\begin{align*}
										\int_{\mathbb{R}^{6}}    	i\hbar p_i  \Theta^-_{K_j}[f] h   \dif p\dif x
										=& \frac{ 1}{(2\pi)^{3}}	\int_{\mathbb{R}^{12}}      \left[p_i'   \delta_- K_j \left(x, \zh \right) - i\frac{\hbar}{2}\dpp{}{x_i}\delta_+ K_j \left(x, \zh \right) \right] f \left(x,p' \right) e^{- i\zh(p -p') }  h(x,p)  \mathrm{d} \zh   \mathrm{d} p' \mathrm{d}p\mathrm{d} x.
									\end{align*}
									Similarly, the second term gives
									\begin{align*}
										\int_{\mathbb{R}^{6}}    - i \frac{\hbar}{2}  \Theta^+_{K_j}\left[\dpp{f}{x_i}\right]   h   \dif p\dif x
										=	&     \frac{ i\hbar}{2(2\pi)^{3}}  	\int_{\mathbb{R}^{12}}     \left[\dpp{}{x_i} \delta_+ K_j \left(x ,\zh  \right) \right] f  \left(x,p' \right) e^{- i\zh(p -p') }  h(x,p)  \dif \zh   \dif p'     \dif p\dif x  	 \\
										&  +    \frac{i \hbar}{2(2\pi)^{3}}  	\int_{\mathbb{R}^{12}}     \left[\delta_+ K_j \left(x ,\zh  \right) \right] f \left(x,p' \right) e^{- i\zh(p -p') }  \dpp{h(x,p)}{x_i}   \dif \zh   \dif p'     \dif p\dif x\;. 
									\end{align*}
									Taking the sum, we obtain
									\begin{align*}
										\int_{\mathbb{R}^{2d}}   	\left[  p_i K , f \right]_{\#} h\dif p\dif x
										=&
										\int_{\mathbb{R}^{2d}}      f \left(x,p \right) \left[ - i \hbar  p_i \Theta^-_K  \left[h \right]  	  + 	i \frac{\hbar}{2}    \Theta^+_K  \left[\dpp{h}{x_i} \right] \right]   \dif p\dif x\;,
									\end{align*}
									which verifies Eq \eqref{Wig_sym_01}. The calculations for Eq. \eqref{Wig_sym_02} proceed similarly. 
									By using the definition of Rashba Hamiltonian in Eq. \eqref{Rash_Ham} and the expansion of the Wigner function on the Pauli basis, we have
									\begin{align*}
										&\left\langle 	[\mathcal{H}^{K} ,f]_\#  ,h \right\rangle= \hbar \sum_{i,j,k=1}^3 \ze_{ijk} \;	\left\langle 	[  p_i K_j  \zs_k, f ]_\#  ,h \right\rangle=\hbar \sum_{i,j,k=1}^3 \ze_{ijk} \;\sum_{r,s=0}^3	\left\langle 	[  p_i K_j  \zs_k,\zs_s f_s ]_\#  ,\zs_r h_r \right\rangle\\
										%
										=&\hbar \sum_{i,j,k=1}^3 \ze_{ijk} \left[ \int_{\mathbb{R}^{2d}}\left( 	 \left[ p_i K_j ,f_0 \right]_\# h_k + \left[p_i K_j ,  f_k \right]_\# h_0  \right)  \dif p\dif x+\sum_{r,s=1}^3	i\ze_{ksr} 	 \int_{\mathbb{R}^{2d}} \left\{ p_i K_j  ,  f_s \right\}_\#  h_r \dif p\dif x     \right]\;.
									\end{align*}
									Using Eqs. \eqref{Wig_sym_01}-\eqref{Wig_sym_02}, we obtain Eq. \eqref{symm_H}.
									\endproof
									
									We study the well-posedness of the adjoint equation, by substituting the final value $f_T$ with $\chi_R \;f_T$, where have introduced the smooth cut-off function $\chi_R \in C_0^{\infty} (\mathbb{R}^3_x\times\mathbb{R}^3_p)$, such that $\chi_R =1$ if $|x|^2+|p|^2<R^2$, and $\textrm{supp}(\chi_R)\subset B_{2R}$, where $B_R$ denotes the ball with radius $R$. 
									\begin{theorem}\label{teor_wellpos_h_dyn}
										Under AS, let $T,R>0$, and $f_T\in C^\infty(\mathbb{R}^3_x\times\mathbb{R}^3_p)$. The solution of the adjoint Wigner adjoint equation $\disp i\hbar \dpp{h}{t} =  \left[ \mathcal{H},h \right]_{\#}$ with Hamiltonian given by Eq. \eqref{tot_Ham} and final condition $h(.,.,T)=f_T\; \chi_{R}$, has a unique strong solution $h$. Moreover, $\|h(.,.,t)\|_{L^2\left( \mathbb{R}^3_x\times\mathbb{R}^3_p,\mathbb{R}^4 \right)}=\|f_T\; \chi_{R}\|_{L^2\left( \mathbb{R}^3_x\times\mathbb{R}^3_p,\mathbb{R}^4 \right)}$, $h\in C\left([0,T], H^1_p\right)$ and for all $t\in[0,T] $, $h $ satisfies
										\begin{align}
											\| h(\cdot,\cdot,t) \|_{H^1_p} \leq&  \| f_T\; \chi_{R} \|_{H^1_p} e^{C \int_{T-t}^T \left(\| K(.,t)\|_{W^{2,\infty}(\mathbb{R}_x^3)^3} +\| U(.,t)\|_{W^{1,\infty}(\mathbb{R}_x^3)}+\| B(.,t)\|_{W^{1,\infty}(\mathbb{R}_x^3)^3} \right)\dif t} \;.\label{H1_bound_sol_h}
										\end{align}
										In the limit $\hbar \rightarrow 0$, $h(t)$ is a solution of the associated adjoint Liouville-Rashba-Zeeman problem
										\begin{align}\label{Adj_Lio_Rashb_prob}
											\left\{\begin{array}{lll}
												\disp	\dpp{h_0}{t}  +  \frac{p}{m} \cdot\nabla_x h_0 -\nabla_xU \cdot\nabla_p h_0 =0 &\textrm{on } \mathbb{R}_x^3\times\mathbb{R}^3_p\times[0 ,T] \\[6pt]
												\disp	\dpp{\vec{h}}{t} +\frac{p}{m} \cdot\nabla_x \vec{h} -\nabla_x U \cdot\nabla_p \vec{h} = 2    ( p\wedge K - B  )\wedge \vec{h} &\textrm{ on } \mathbb{R}_x^3\times\mathbb{R}_p^3\times[0 ,T] \\[6pt]
												\left.\vec{h}_i \right|_{t=T}=  \frac{\chi_{R}}{2}\mathrm{tr} \left\{\zs_i f_T \right\} &\textrm{ on } \mathbb{R}_x^3\times\mathbb{R}_p^3\\
												\left.h_0 \right|_{t=T} = \frac{\chi_{R}}{2}\mathrm{tr} \left\{ f_T\right\}  &\textrm{ on } \mathbb{R}_x^3\times\mathbb{R}_p^3
											\end{array}
											\right.
										\end{align}
										where $h_0:\mathbb{R}_x^3\times\mathbb{R}_p^3\rightarrow  \mathbb{R}$ and $\vec{h}=\mathbb{R}_x^3\times\mathbb{R}_p^3\rightarrow  \mathbb{R}^3$.  
									\end{theorem}
									\proof Since the evolution equation of the adjoint and the Wigner dynamics coincide, the theorem follows from straightforward modifications of the proof of Theorem \ref{teor_wellpos_W_dyn}. In particular, the cut-off ensures that the final datum belongs to $H^1_p$. The conservation of the L2 norm follows directly from the anti-symmetry of the Moyal commutator of the Hamiltoinan expressed by Lemma \ref{lemma_anti_symm_H}. Equation \eqref{H1_bound_sol_h} mimics Eq. \eqref{H1_bound_sol_wig}. Concerning the classical limit, Cauchy-Schwartz estimates $\left|\left\langle h^\ze ,\mathsf{h} \right\rangle_{\mathrm{HS}} \right| \leq  \| \chi_R\; f_T  \|_{L^2(\mathbb{R}_x^3\times \mathbb{R}_p^3)}  \| \| \mathsf{h} \|_{L^2(\mathbb{R}_x^3\times \mathbb{R}_p^3)}  $ 
									with $\mathsf{h}$ Schwartz function, let us take the limit. Calculations similar to what was done in the proof of Th. \eqref{Theor_lim_eq_Wig_f} guarantee that the limit function solves the Liouville-Rashba-Zeeman problem of Eq. \eqref{Adj_Lio_Rashb_prob}.
									\endproof 
									
									\subsection*{$III:$ Concentration limit and convergence of the control equation}
									The aim of this section is to apply the results obtained in the previous two steps concerning the limit of the Wigner and adjoint functions, in the integral term appearing in the last equation of the OC problem \eqref{Opt_cond_Q} providing the parameters associated with the optimal control. As a final result, we write the control function as a solution of an ODE problem, which can be directly compared with the analogous equation of the OC (Eq. \eqref{char_adj_1D_03}) for the deterministic case.
									
									\begin{lemma}\label{lemma_lim_comp_cont_02}
										Let $(\vec{h},h_0)$, where $\vec{h}:\mathbb{R}^3_x\times\mathbb{R}^3_p\times\mathbb{R}_t \rightarrow \mathbb{R}^3$, $h_0:\mathbb{R}^3_x\times\mathbb{R}^3_p\times\mathbb{R}_t \rightarrow \mathbb{R}$ be a solution of the classical Liouville-Rashba-Zeeman equation 
										\begin{align}
											&\dpp{\vec{h}}{t}+\frac{p}{m}\cdot \nabla_x \vec{h} + E(x,t)\cdot \nabla_p \vec{h}+  B(x,p) \wedge \vec{h} =0& \textrm{ on } \mathbb{R}^3_x\times\mathbb{R}
											_p^3\times[0 ,T
											]\label{Liou_h_traj_01}
											\\
											&\dpp{h_0}{t}+\frac{p}{m}\cdot \nabla_x h_0 + E(x,t)\cdot \nabla_p h_0 =0& \textrm{ on } \mathbb{R}^3_x\times\mathbb{R}
											_p^3\times[0 ,T
											]\label{Liou_h0_traj_01}
											\\
											&\left.(\vec{h}, h_0)\right|_{t=T}=(\vec{f}_T,[f_T]_0)  & \textrm{ on } \mathbb{R}^3_x\times\mathbb{R}
											_p^3
										\end{align}
										with $f_T \in C^0\left(\mathbb{R}_x\times\mathbb{R}
										_p\right)^4$, the total magnetic field is the sum of the external field and the Rashba term $B(x,p)=2(B(x)-K(x)\wedge p)$, the total electric field $E(x,t)=-\nabla_x U (x,u(t))$, where the dependence on time arises form the control parameters $u$ which here are assumed to be given, and AS holds true. 
										We have
										\begin{align}
											&d(t)\cdot  \left.	\dpp{\vec{h}}{p_i}\right|_{ (x(t),p(t))}+\left.	\dpp{h_0}{p_i}\right|_{ (x(t),p(t))}    = x^h_i( t) & t \in [0,T],
										\end{align}
										where $ (x(t),p(t))$ denotes the phase-space flow given by Eqs. \eqref{Liou_h_traj_01}-\eqref{Liou_h0_traj_01}, with final condition $	(x(T),p(T) )= (\overline{x},\overline{p}) $, the vector $d$ rotates around the total magnetic field evaluated along the classical trajectory
										\begin{align}
											\left\{
											\begin{array}{ll}
												\disp 	\dpt{\; d}{t}= B\left(x(t),p(t)\right) \wedge d  & \textrm{on } t\in [0,T] \\[6pt]
												\disp	 d(T)=\overline{d} & \textrm{in } t=T
											\end{array}
											\right.\label{cauchy_d_T}
										\end{align}
										and $x^h:[0,T]\rightarrow \mathbb{R}^3$ is obtained from  the Cauchy problem  
										\begin{align}
											\left\{
											\begin{array}{ll}
												\disp \dpt{x^h_i }{t} =   - \frac{p^h_i}{m}   +  \dpp{B}{p_i} \cdot\zh^h
												&i=1,2,3,  \textrm{ on } [0 ,T] \\[6pt]
												\disp \dpt{p^h_i}{t} =     - \dpp{E}{x_i}\; x^h - \dpp{ B }{x_i}  \cdot\zh^h  &i=1,2,3,  \textrm{ on } [0 ,T] \\[6pt]
												\disp \dpt{\zh^h}{t} =   - B \wedge \zh^h  &\textrm{ on } [0 ,T] \\[6pt]
												\disp x^h_i(T) =d\cdot  \left. \dpp{\vec{f_T}}{p_i}\right|_{(x,p)=(\overline{x},\overline{p})} +\left.\dpp{[f_T]_0}{p_i}\right|_{(x,p)=(\overline{x},\overline{p})}& \textrm{ in } t=T\\[8pt]
												\disp p_i^h(T) =d\cdot  \left. \dpp{\vec{f_T}}{x_i}\right|_{(\overline{x},\overline{p})}  + \left. \dpp{[f_T]_0}{x_i}\right|_{(\overline{x},\overline{p})} & \textrm{ in } t=T\\[6pt]
												\zh^h(T) = \vec{f_T}(\overline{x},\overline{p}) \wedge d  & \textrm{ in } t=T
											\end{array}\right.\label{cahr_sys_lim_h0}
										\end{align}
									\end{lemma}
									\proof 
									Equations \eqref{Liou_h_traj_01}-\eqref{Liou_h0_traj_01} have the same phase-space flow $(x(t),p(t))$ associated to the characteristic system 
									$ (x(t),p(t)):\left\{
									\begin{array}{l}
										\dpt{x}{t}=	\frac{p}{m} \\
										\dpt{p}{t}=	E(x(t),t)  
									\end{array}
									\right.$, where $E=-\nabla_x U$ which, for the assumption on $U$, admits a global solution. We determine the solution of Eq. \eqref{Liou_h_traj_01} using the characteristic method.
									It is convenient to derive the characteristic equations in a more general case. We consider the following system of ODE 
									\begin{align}
										\dpp{h}{t} +\sum_{i=1}^{6} \dpp{h}{X^i} g_i(X)+ H (h,X)=0 \;.\label{char_meth_01}
									\end{align}
									Here, we have defined the collective phase-space coordinates $X\doteq (x,p)$, i. e. $X_i = \left\{ \begin{array}{ll}
										x_i & \textrm{for } i=1,2,3 \\
										p_i &  \textrm{for }i=4,5,6
									\end{array}\right.$. Roman indices refer to spatial and momentum variables, Greek indices $\mu=1,2,3$ refer to spin degrees of freedom. Equations \eqref{Liou_h_traj_01} have the form of  Eq. \eqref{char_meth_01} when $g_i = \left\{ \begin{array}{ll}
										\frac{p_i}{m} & \textrm{for } i=1\leq i\leq3  \\
										E_i &  \textrm{for } i=4\leq i\leq6
									\end{array}\right.$, and
									$H(z,x,p) = B \wedge z$. The characteristic equations are (the details are given in Appendix \ref{App_char_sys})
									\begin{align}
										\left\{
										\begin{array}{ll}
											\disp	\dpt{X_j}{s}=	g_j&  \\[6pt]
											\disp	\dpt{P_{X_j} }{s}  = -\sum_{i=1}^{6} P_{X_i}  \dpp{g_i}{X_j} - \dpp{H}{X_j} - \sum_{\nu=1}^3 \dpp{H}{z^{\nu}}P^{\nu}_{X_j} & \\[6pt]
											\disp 	\dpt{z}{s}=- H (z,X)
										\end{array}
										\right.
										\label{char_meth_02}
									\end{align}
									with $j=1,\ldots,6$, we have defined $P_{x_i} \doteq \dpp{h}{x_i}$, $P_{p_i} \doteq \dpp{h}{p_i}$, and $z= h (x(t),p(t),t)$ is the solution evaluated along the characteristic set. The first and third lines of Eq. \eqref{char_meth_02} provide the characteristic equations for the trajectories and for the solution, respectively
									\begin{align*}
										\left\{
										\begin{array}{ll}
											\disp 	\dpt{x}{t}=	\frac{p}{m} & \textrm{on } t\in [0,T]\\[6pt]
											\disp	\dpt{p}{t}=	E & \textrm{on } t\in [0,T]\\[6pt]
											\disp 	\dpt{z}{t}=  - B \wedge z  & \textrm{on } t\in [0,T] \\[6pt]
											\disp	(x(T),p(T) )= (\overline{x},\overline{p}), \; z(T)=f_T (\overline{x},\overline{p})  & \textrm{in } t=T
										\end{array}
										\right.\;.
									\end{align*}
									The second line of Eq. \eqref{char_meth_02} gives the equations for the adjoint variables
									\begin{align*}
										\left\{
										\begin{array}{ll}
											\dpt{P_{x_i}}{s}  = -\sum_{j=1}^3 P_{p_j}   \dpp{E_j}{x_i}  -B\wedge P_{x_i} -\dpp{ B }{x_i} \wedge z   & \textrm{on } t\in [0,T]\\[4pt]
											\dpt{P_{p_i}}{s}  = - \frac{P_{x_i} }{m}   -\dpp{ B }{p_i}\wedge z  -B \wedge P_{p_i}  & \textrm{on } t\in [0,T]\\[4pt]
											P_{p_i}(T) =  \left. \dpp{f_T}{p_i}\right|_{(\overline{x},\overline{p})}, \; P_{x_i}(T) =\left. \dpp{f_T}{x_i}\right|_{(\overline{x},\overline{p})} \;  & \textrm{ in } t=T
										\end{array}
										\right.\;,
									\end{align*}
									with $i=1,2,3$. 
									We define ${x}^v_i\doteq d(t)\cdot  \left.	\dpp{h}{p_i}\right|_{ (x(t),p(t))} = (P_p\cdot d)$, $p^v_i = (P_{x_i}\cdot d)$ and $\zh = z\wedge d$, where $d$ is the solution of the Cauchy problem \eqref{cauchy_d_T} with prescribed final value. We obtain the set of equations
									\begin{align*}
										\left\{
										\begin{array}{ll}
											\disp \dpt{x^v_i }{t} =   - \frac{p^v_i}{m}   +  \dpp{B}{p_i} \cdot\zh
											&i=1,2,3,  \textrm{ on } [0 ,T] \\[6pt]
											\disp \dpt{p^v_i}{t} =     - \dpp{E}{x_i}\; x^v - \dpp{ B }{x_i}  \cdot\zh  &i=1,2,3,  \textrm{ on } [0 ,T] \\[6pt]
											\disp \dpt{\zh }{t} =   - B \wedge \zh  &\textrm{ on } [0 ,T] \\[6pt]
											\disp x^v_i(T) =d\cdot  \left. \dpp{f_T}{p_i}\right|_{(x,p)=(\overline{x},\overline{p})} & \textrm{ in } t=T\\[6pt]
											\disp p_i^v(T) =d\cdot  \left. \dpp{f_T}{x_i}\right|_{(\overline{x},\overline{p})}  & \textrm{ in } t=T\\[6pt]
											\zh(T) = f_T(\overline{x},\overline{p}) \wedge d & \textrm{ in } t=T
										\end{array}\right.
									\end{align*}
									The same arguments apply to Eq. \eqref{Liou_h0_traj_01}, where we define ${x}^0_i\doteq  \left.	\dpp{h_0}{p_i}\right|_{ (x(t),p(t))} = P_{p_i}$ and $p^0_i = P_{x_i}$. We obtain
									\begin{align*}
										\left\{
										\begin{array}{ll}
											\disp \dpt{x^0_i }{t} =   - \frac{p^0_i}{m}   
											&i=1,2,3,  \textrm{ on } [0 ,T] \\[6pt]
											\disp \dpt{p^0_i}{t} =     - \dpp{E}{x_i}\; x^0  &i=1,2,3,  \textrm{ on } [0 ,T] \\[6pt]
											\disp x^0_i(T) =\left. \dpp{[f_T]_0}{p_i}\right|_{(x,p)=(\overline{x},\overline{p})} & \textrm{ in } t=T\\[6pt]
											\disp p_i^0(T) =  \left. \dpp{[f_T]_0}{x_i}\right|_{(\overline{x},\overline{p})}  & \textrm{ in } t=T
										\end{array}\right.
									\end{align*}
									Taking $x^h=x^v+x^0$ and $p^h=p^v+p^0$ we have the thesis. 
									\endproof
									As a last step for the proof of Th. \ref{Theor_lim_OC}, we collect the previous results and remove the cut-off on the final value of the adjoint variable 
									\begin{corollary}\label{Corol_lim_cont}
										For any $(\overline{x},\overline{p})\in \mathbb{R}_x^3\times \mathbb{R}_p^3$, it is always possible to find a cut-off $R>0$ such that 
										\begin{align*} 
											\lim_{\hbar \rightarrow 0}  \frac{1}{2} \textrm{tr} \int_{\mathbb{R}_x^3\times\mathbb{R}_p^3}   \Theta^-_{\dpp{U}{u_j}} [h] \;f    \dif x  \dif p =x^h_j\;,
										\end{align*}
										where $j=1,2,3$, and $x^h$ is obtained as a solution of the characteristic system given by  Eq. \eqref{cahr_sys_lim_h0}.
									\end{corollary}
									\proof 
									For any $R>0$, from Theorem \ref{Theor_lim_eq_Wig_f} we have that when $\hbar$ goes to zero, the Wigner function $f$ concentrates to a Dirac's delta localized along the classical trajectories \eqref{Char_traj_1D}, while according to Theorem \ref{teor_wellpos_h_dyn}, the adjoint function $h$ belongs to $H^1_p$. 
									Moreover, estimations in Theorem \ref{Theor_lim_eq_Wig_f} ensure that the pseudo-differential operator $\Theta^-$ can be approximated by a derivative operator. Thus, we have
									\begin{align} 
										\lim_{\hbar \rightarrow 0}  \frac{1}{2} \textrm{tr} \int_{\mathbb{R}_x^3\times\mathbb{R}_p^3}   \Theta^-_{\dpp{U}{u_j}} [h] \;f    \dif x  \dif p =& \frac{1}{2} \textrm{tr} \int_{\mathbb{R}_x^3\times\mathbb{R}_p^3} \dpp{E}{u_j}\cdot  \nabla_p h \;f    \dif x  \dif p \nn \\ 
										= &\sum_{i} \dpp{E_i}{u_j} \left.\dpp{\vec{h}}{p_i}\right|_{(x(t),p(t))}  d(t) +\sum_{i} \dpp{E_i}{u_j} \left.\dpp{h_0}{p_i}\right|_{(x(t),p(t))}\label{lim_cont_eq}
									\end{align}
									where, the dipole $d$ precesses according to Eq. \eqref{Char_traj_1D}. By the conservation of the energy, we have that the momentum is bounded $|p|^2 \leq {\overline{p}}^2 +4m \| U\|_{L^\infty}$, and the position
									$|x|\leq |\overline{x}| + T \sqrt{ \left(\frac{\overline{p}}{m}\right)^2+ \frac{4  \| U\|_{L^\infty}}{m}}$. Due to the fact that the characteristic set is bounded, it is always possible to find $R$ sufficiently large that the classical trajectories are contained inside the ball of radius $2R$. In this case, Lemma \ref{lemma_lim_comp_cont_02} ensures that the right side of Eq. \eqref{lim_cont_eq} is equal to the component $x^h$ of the solution of the characteristic system \eqref{cahr_sys_lim_h0}.
									\endproof
									\section{Conclusions}
									
									The classical limit of an optimal control problem for the quantum dynamics of a charged particle gas with spin has been investigated. The problem is studied in an unbounded domain in the presence of external magnetic and potential fields and Rashba interaction. The dynamics of the system is controlled by including in the model a set of parameters modulating the total potential and which are designed with the aim to miminize the second momenta of the position and velocity of the particle distribution and the spin projection along a given direction, at the final time. A cost function representing the amount of external energy required to control the system and should be minimized is also considered. As a main result, we have shown that in the classical limit where the Planck constant goes to zero, the problem of optimal control of the quantum dynamics converges to the analogous optimal control problem for the classical deterministic case. 
									
									\section{Acknowledgments} 
									The work was developed under the auspices of
									GNFM (INdAM). This project received funding from MUR - Next Generation EU project, PRIN 2022 and TRANFORM project, PNR 2021-2027.
									
									\appendix
									\section{Derivation of Eq. \eqref{Wig_f0_fk_01}}\label{App_der_eq_wig}
									We derive the Wigner Eq. \eqref{Wig_f0_fk_01}, which describes the dynamics of the quantum system. 
									The trace of the Wigner evolution equation in the Moyal representation $i\hbar \dpp{f}{t}=[\mathcal{H},f]_\#$, gives 
									$i\hbar \dpp{f_k}{t} =\sum_{i,j=0}^3 \frac{1}{2} \textrm{tr}\left[\zs_k\zs_i\zs_j \right]  \left(\mathcal{H}_i\# f_j- f_i\#   \mathcal{H}_j\right)$, 
									where we have decomposed the Wigner function and the Hamiltonian symbol on the Pauli basis, $f=\sum_{i=0}^3 f_i \zs_i$, $\mathcal{H}=\sum_{i=0}^3 \mathcal{H}_i \zs_i$ with  $f_i=\frac{1}{2}\textrm{tr}\left\{f\zs_j \right\} $, $\mathcal{H}_i=\frac{1}{2}\textrm{tr}\left\{\mathcal{H} \zs_j \right\} $. We obtain 
									\begin{align*}
										\hbar \dpp{f_0}{t} =&\sum_{i=0}^3\textrm{Im}\left[\mathcal{H}_i, f_i\right]_\#   \\
										\hbar \dpp{f_k}{t} =& \textrm{Im}\left[\mathcal{H}_0, f_k\right]_{\#}+ \textrm{Im}\left[\mathcal{H}_k, f_0\right]_{\#}+\sum_{i,j=1}^3  \ze_{ijk}\textrm{Re}\left\{\mathcal{H}_i, f_j\right\}_{\#} & k=1,2,3.
									\end{align*}
									We detail the calculations of the Rashba Hamiltonian $
									\mathcal{H}^{K} = \hbar\left( p \wedge K\right) \cdot \zs$. 
									In order to simplify the notation, in the next calculations we drop the index notation of $f$ and $K$ and assume that $f$ and $K$ represent any component $f_i$ and $K_j$, respectively. From the definition of the Moyal product 
									\begin{align}
										F \# G \left( x, p\right)=&    \frac{ 1}{(2\pi)^{6}}  \int_{\mathbb{R}^{12}} F\left(x -\frac{\hbar \zh'}{2},p+ \frac{\hbar\zm'}{2} \right)  G \left(x',p' \right) e^{i\mu'(x-  x')+i\zh'(p -p') } \dif\mu' \dif \zh'   \dif x'\dif p'\;,\label{Moyal_prod}
									\end{align}
									we have 
									\begin{align}
										\left[  p_i K , f \right]_{\#} =& 	  p_i K  \# f -	   f\#  p_i K   =
										i\hbar p_i  \Theta^-_{K}[f]   - i \frac{\hbar}{2}  \Theta^+_{K}\left[\dpp{f}{x_i}\right]  \label{Moy_com_pK}	\\
										\left\{ p_i K  ,  f \right\}_{\#} =&  	p_i K  \# f  +	  f\#  p_i K   =	 p_i  \Theta^+_{K}[f]   +  \frac{\hbar^2}{2}  \Theta^-_{K}\left[\dpp{f}{x_i}\right] \;, \label{Moy_acom_pK}
									\end{align}
									where the pseudo-differential operators $\Theta^\pm$ are defined in Eq. \eqref{Theta_pm}.   
									We have
									\begin{align*}
										\frac{1}{\hbar}\sum_{i=0}^3\textrm{Im}\left[\mathcal{H}_i^{K}, f_i\right]_\# 
										=&\sum_{i,j,k=1}^3  \ze_{ijk}\left(  \hbar  p_i 	\Theta^-_{K_j}[f_k]     -  \frac{\hbar}{2}  \Theta^+_{K_j}\left[\dpp{f_k}{x_i}\right] \right)\\
										\frac{1}{\hbar}\textrm{Im}\left[\mathcal{H}_k^K, f_0\right]_{\#} 
										=&\sum_{i,j=1}^3  \ze_{ijk}\left(  \hbar   p_i 	\Theta^-_{K_j}[f_0]     -  \frac{\hbar}{2}  \Theta^+_{K_j}\left[\dpp{f_0}{x_i}\right] \right)\\
										\frac{1}{\hbar}\sum_{i,j=1}^3  \ze_{ijk}\textrm{Re}\left\{\mathcal{H}_i, f_j\right\}_{\#} 
										=& \sum_{i,j,l,s=1}^3 \ze_{ils} \;    \ze_{ijk} \left(    p_l 	\Theta^+_{K_s}[f_j]    + \frac{ \hbar^2}{2}  \Theta^-_{K_s}\left[\dpp{f_j}{x_l}\right] \right)\;.
									\end{align*}
									In the same way, the Zeeman term $-\hbar B(x) \cdot \zs$ gives $ \sum_{i=0}^3\textrm{Im}\left[ B_i , f_i\right]_\# = \sum_{i=1}^3 \hbar \Theta^-_{B_i}[f_i] $, $ 
									\textrm{Im}\left[B_k, f_0\right]_{\#}= \hbar \Theta^-_{B_k}[f_0]  $, and $\sum_{i,j=1}^3  \ze_{ijk}\textrm{Re}\left\{B_i, f_j\right\}_{\#}  = \sum_{i,j=1}^3 \ze_{ijk}\Theta^+_{B_i}[f_j]$, with $k=1,2,3 $.
									Using the previous equations, completed with the kinetic and potential terms for which the calculation is standard, we obtain Eq. \eqref{Wig_f0_fk_01}. 

							\section{Weak form of the Wigner dynamics}\label{app_der_w_form}
							We derive the Wigner equation in the weak formulation. We multiply Eq. \eqref{Wig_f0_fk_01} by a test function $h$  and we integrate in time and in phase-space coordinates. We obtain
							\begin{align*}
								\disp	\frac{1}{2}\mathrm{tr} \int_{\mathbb{R}^{6}\times \mathbb{R}^{+}}\left(\dpp{f}{t} h + \left(\frac{p}{m} \cdot\nabla_xf \right) h -\Theta^-_U [f] h\right) \dif p \dif x  \dif t =& I+II+III\; ,
							\end{align*}
							where we have defined
							\begin{align*}
								I\doteq & \sum_{k=1}^3 \int_{\mathbb{R}^{6}\times \mathbb{R}^{+}}\left(A^+[f_k]_k h_0  +A^+[f_0]_k h_k  \right)\dif p \dif x  \dif t \\
								II\doteq & \sum_{i,j,k=1}^3   \ze_{ijk} \int_{\mathbb{R}^{6}\times \mathbb{R}^{+}} A^-[f_j]_i h_k  \dif p \dif x  \dif t \\
								III\doteq &-\sum_{k=1}^3 \int_{\mathbb{R}^{6}\times \mathbb{R}^{+}}   \left(   \hbar\Theta^-_{ B_k}[f_k] h_0   +   \hbar\Theta^-_{ B_k}[f_0]h_k +  \sum_{i,j=1}^3 \ze_{ijk}\, \Theta^+_{ B_i}[f_j] h_k\right)\dif p \dif x  \dif t \;.
							\end{align*}
							Using $ \int_{\mathbb{R}^{d}} \Theta^-_{ U}[f_i] h_j \dif p=-\int_{\mathbb{R}^{d}} f_i\Theta^-_{ U}[h_j] \dif p$, and $ \int_{\mathbb{R}^{d}} \Theta^+_{ U}[f_i] h_j \dif p=\int_{\mathbb{R}^{d}} f_i\Theta^+_{ U}[h_j] \dif p$, we obtain
							\begin{align*}
								I&= - \frac{ \hbar}{2}     \sum_{i,j,k=1}^3  \ze_{ijk} \mathrm{Re }\; \mathrm{tr}\int_{\mathbb{R}^{6}\times \mathbb{R}^{+}}    \zs_k   \left(	\Theta^-_{K_j}[ p_i h]+ \frac{1}{2}   	\Theta^+_{K_j}\left[\dpp{ h}{x_i}\right]        \right)f\dif p \dif x  \dif t \\
								II&=   -\frac{1}{2} \sum_{k,i,j=1}^3  \ze_{ijk}\;   \mathrm{Im }\;  \mathrm{ tr}  \int_{\mathbb{R}^{6}\times \mathbb{R}^{+}}     \zs_k   \left(    \Theta^+_{K_{j}}[p_{i} h ] +\frac{\hbar^2}{2}   \dpp{}{x_{i}}\Theta^-_{K_{j}}\left[h_i \right] \right) f  \dif p \dif x  \dif t \\
								III&=	\frac{1}{2} \int_{\mathbb{R}^{6}\times \mathbb{R}^{+}}\left(\hbar\mathrm{Re }\;  \mathrm{tr} \left\{ f\Theta^-_{ B}[ h]  \right\}  +\mathrm{Im }\; \mathrm{tr}   \left\{    f  \Theta^+_{ B}[h]  \right\} \right)\dif p \dif x  \dif t \;.
							\end{align*}
							This gives Eq. \eqref{Wig_f0_w_01}.  

							\section{Method of the characteristics}\label{App_char_sys}
							We consider the following set of equations
							\begin{align}
								\sum_{i=0}^{2d} \dpp{h}{X^i} g_i(X)+ H (z,X)=0 \;.\label{char_meth_01}
							\end{align}
							Where $h:\mathbb{R}^d_x\times \mathbb{R}^d_p\rightarrow \mathbb{R}^3$ we have defined the collective phase-space coordinates $X\doteq (t,x,p)$ with $X_0=t$, $X_i=x_i$, $i=1,\ldots,d$ and $X_i=p_i$ $i=d+1,\ldots,2d$. Roman indices $i=1,d$, refer to spatial momentum variables, Greek indices $\mu=1,\ldots ,3$ refer to spin degrees of freedom.
							Using standard notation, we write the partial diff equations as
							\begin{align*}
								\sum_{i=0}^{2d} P_{X_i}  g_i + H (z,X)=0\;,
							\end{align*}
							where $P_t \doteq \dpp{h}{t}$, $P_{x_i} \doteq \dpp{h}{x_i}$, $P_{p_i} \doteq \dpp{h}{p_i}$ and $z= h (x(t),p(t),t)$ is the solution evaluated along the characteristic set. We take the total derivative with respect to the arc parameter $s=t$ of Eq. \eqref{char_meth_01}
							\begin{align*}
								\sum_{i=1}^{2d} \dpt{P_{X_i}}{s} g_i  +\sum_{i,j}^{2d}  P_{X_i}  \dpp{g_i}{X_j} \dpt{X_j}{s}+ \dpp{H^{\mu}}{X_j} \dpt{X_j}{s}+ \sum_{\nu}\dpp{H^{\mu}}{z^{\nu}}\dpt{z^{\nu}}{s}=0 
							\end{align*}
							and we use $
							\dpt{z^{\nu}}{s}= \sum_s \dpp{f^{\nu}}{X_s} \dpt{X_s}{s} =\sum_s P^{\nu}_{X_s} \dpt{X_s}{s} $. We obtain
							\begin{align*}
								\sum_{j=0}^{2d}\left(\dpt{P_{X_j} }{s}  +\sum_{i=0}^{2d} P_{X_i}  \dpp{g_i}{X_j} + \dpp{H}{X_j} + \sum_{\nu=1}^3\dpp{H}{z^{\nu}}P^{\nu}_j\right)\dpt{X_j}{s}=0 \;,
							\end{align*}
							where we have assumed $g_j=\dpt{X_j}{s}$. We obtain the characteristic set of equations
							\begin{align*}
								\left\{
								\begin{array}{ll}
									\disp	\dpt{X_j}{s}=	g_j& i=1,\ldots,2d\\[6pt]
									\disp	\dpt{P_{X_j} }{s}  = -\sum_{i=0}^{2d} P_{X_i}  \dpp{g_i}{X_j} - \dpp{H}{X_j} - \sum_{\nu=1}^3 \dpp{H}{z^{\nu}}P^{\nu}_{X_j} & i=1,\ldots,2d\\[6pt]
									\disp 	\dpt{z}{s}=- H (z,X)
								\end{array}
								\right.\;.
							\end{align*}

							\section{Control on the H1 norm of the Wigner function} \label{App_est_d_xmu}
							We provide a few expressions that are applied in Theorem \ref{teor_wellpos_W_dyn} to the control of the H1 norm of the Wigner function. By using \eqref{moy_prod_pKf_01} we obtain 
							\begin{align*}
								&\left[  p_i K_j , \dsp{}{x_\mu} f_0 \right]_{\#}	-\dsp{}{x_\mu}	\left[  p_i K_j , f_0 \right]_{\#} = \\
								&	i\hbar p_i  \Theta^-_{K_j}\left[\dsp{}{x_\mu}f_0\right]
								-i\hbar\dsp{}{x_\mu}   p_i \Theta^-_{K_j}[f_0]   + i \frac{\hbar}{2} \dsp{}{x_\mu}\Theta^+_{K_j}\left[\dpp{f_0}{x_i}\right]    - i \frac{\hbar}{2}  \Theta^+_{K_j}\left[\dpp{}{x_i}\dsp{}{x_\mu} f_0\right] \\
								=&
								-i\hbar p_i \Theta^-_{\dsp{K_j}{x_\mu}  }[f_0]  -2i  \hbar   p_i \Theta^-_{\dpp{K_j}{x_\mu}  }\left[\dpp{f_0}{x_\mu} \right]  + i \frac{\hbar}{2} \Theta^+_{\dsp{}{x_\mu}K_j}\left[\dpp{f_0}{x_i}\right]     + i  \hbar \Theta^+_{\dpp{K_j}{x_\mu}}\left[\dpp{}{x_\mu}\dpp{f_0}{x_i}\right]\;.
							\end{align*}
							We estimate
							\begin{align*}
								&\left| \int_{\mathbb{R}^{2d}} \hbar p_i \Theta^-_{\dsp{K_j}{x_\mu}  }[f_0] f_k\dif x\dif p\right|  \leq  2 \left\| \dsp{K}{x_\mu}  \right\|_{L^\infty (\mathbb{R}^{d}_{x})} 	 \left\| p_i 	 f_k    \right\|_{L^2 (\mathbb{R}^{2d}_{x,p})}  \left\|   f_0 \right\|_{L^2 (\mathbb{R}^{2d}_{x,p})}  \\
								&\left| \int_{\mathbb{R}^{2d}} 2  \hbar   p_i \Theta^-_{\dpp{K_j}{x_\mu}  }[\dpp{f_0}{x_\mu} ] f_k\dif x\dif p\right|  \leq  4 \left\| \dpp{K}{x_\mu}  \right\|_{L^\infty (\mathbb{R}^{d}_{x})} 	 \left\| p_i 	 f_k    \right\|_{L^2 (\mathbb{R}^{2d}_{x,p})}  \left\|   \dpp{f_0}{x_\mu}  \right\|_{L^2 (\mathbb{R}^{2d}_{x,p})}  \\
								&\left| \int_{\mathbb{R}^{2d}}  \frac{\hbar}{2} \Theta^+_{\dsp{}{x_\mu}K_j}\left[\dpp{f_0}{x_i}\right]    f_k\dif x\dif p\right|  \leq \hbar \left\| \dsp{K}{x_\mu}  \right\|_{L^\infty (\mathbb{R}^{d}_{x})} 	 \left\|	 f_k    \right\|_{L^2 (\mathbb{R}^{2d}_{x,p})}  \left\|   \dpp{f_0}{x_\mu}  \right\|_{L^2 (\mathbb{R}^{2d}_{x,p})}  \\
								& \int_{\mathbb{R}^{2d}}    \hbar \Theta^+_{\dpp{K_j}{x_\mu}}\left[\dpp{}{x_\mu}\dpp{f_0}{x_i}\right]     f_k\dif x\dif p =  -\int_{\mathbb{R}^{2d}}    \hbar \Theta^+_{\dsp{K_j}{x_\mu}}\left[\dpp{f_0}{x_i}\right]     f_k\dif x\dif p - \int_{\mathbb{R}^{2d}}    \hbar \Theta^+_{\dpp{K_j}{x_\mu}}\left[\dpp{f_0}{x_i}\right]  \dpp{  f_k}{x_\mu} \dif x\dif p \\
								& \left|\int_{\mathbb{R}^{2d}}    \hbar \Theta^+_{\dpp{K_j}{x_\mu}}\left[\dpp{}{x_\mu}\dpp{f_0}{x_i}\right]     f_k\dif x\dif p \right|\leq 2\hbar\left( \left\| \dsp{K}{x_\mu}  \right\|_{L^\infty (\mathbb{R}^{d}_{x})}  + \left\| \dpp{K}{x_\mu}  \right\|_{L^\infty (\mathbb{R}^{d}_{x})} 	 \right)\times	\\
								&\left( \left\|	 f_k    \right\|_{L^2 (\mathbb{R}^{2d}_{x,p})}  + \left\|	 \dpp{f_k}{x_\mu}    \right\|_{L^2 (\mathbb{R}^{2d}_{x,p})}  \right)  \left\|   \dpp{f_0}{x_i}  \right\|_{L^2 (\mathbb{R}^{2d}_{x,p})}.
							\end{align*}
							The anticommutator terms can be estimated in the same way and we obtain
							\begin{align*}
								&\dpt{}{t}\left\langle f,  -\dsp{f}{x_\mu}   \right\rangle  \leq C   \max_\mu  \left(  \left\| \dsp{K}{x_\mu}  \right\|_{L^\infty (\mathbb{R}^{d}_{x})}  + \left\| \dpp{K}{x_\mu}  \right\|_{L^\infty (\mathbb{R}^{d}_{x})}  \right)  \left\|    f   \right\|_{H^1_p } ^2\;.
							\end{align*}
							


\begin{thebibliography}{}
	
	
	
	\bibitem{Nielsen_10} M. A. Nielsen and I. L. Chuang, Quantum Computation and Quantum Information, 2nd ed. (Cambridge University Press,New York, 2010).
	
	\bibitem{Koch_22} C. P. Koch, U. Boscain, T. Calarco, G. Dirr, S. Filipp, S. J. Glaser, R. Kosloff, S. Montangero, T. Schulte-Herbrueggen, D. Sugny, and F. Wilhelm, EPJ Quantum Technol. 9, 19 (2022).
	
	\bibitem{Saffman_19}  M. Saffman, Natl. Sci. Rev. 6, 24 (2019).
	
	\bibitem{Henriet_20} L. Henriet, Quantum 4, 327 (2020).
	
	\bibitem{Bluvstein_23}  D. Bluvstein, S. J. Evered, A. A. Geim, S. H. Li, H. Zhou, T. Manovitz, S. Ebadi, M. Cain, M. Kalinowski, D. Hangleiter, J. P. Bonilla Ataides, N. Maskara, I. Cong, X. Gao, P. S. Rodriguez, T. Karolyshyn, G. Semeghini, M. J. Gullans, M. Greiner, V. Vuleti´ c, and M. D. Lukin, Nature (London) 626, 58 (2023).
	
	\bibitem{Couvert_08}  A. Couvert, T. Kawalec, G. Reinaudi, and D. Gu\'ery-Odelin, Europhysics Letters 83, 13001 (2008).
	
	\bibitem{Murphy_09} M. Murphy, L. Jiang, N. Khaneja, and T. Calarco, Phys. Rev. A 79, 020301(R) (2009).
	
	\bibitem{Chen_11}  X. Chen, E. Torrontegui, D. Stefanatos, J.-S. Li, and J. G. Muga, Phys. Rev. A 84, 043415  (2011).
	
	\bibitem{Cicali_24} C. Cicali, M. Calzavara, E. Cuestas, T. C. R. Zeier, and F. Motzoi, N  arXiv:2412.15173 [quant-ph](2024). 
	
	\bibitem{Rosi_13} S. Rosi, A. Bernard, N. Fabbri, L. Fallani, C. Fort, M. Inguscio, T. Calarco, and S. Montangero, Phys. Rev. A 88, 021601(R) (2013)
	
	
	\bibitem{Torrontegui_11} E. Torrontegui, S. Ibáñez, X. Chen, A. Ruschhaupt, D. Guéry-Odelin, and J. G. Muga, Phys. Rev. A 83, 013415 (2011).
	
	\bibitem{Negretti_13} A. Negretti, A. Benseny, J. Mompart, and T. Calarco, Quantum Inf. Process. 12, 1439 (2013).
	
	\bibitem{GueryOdelin_19} D. Guéry-Odelin, A. Ruschhaupt, A. Kiely, E. Torrontegui, S. Martínez-Garaot, and J. G. Muga, Rev. Mod. Phys. 91, 045001 (2019).
	
	\bibitem{Jaeger_14} G. Jäger, D. M. Reich, M. H. Goerz, C. P. Koch, and U. Hohenester, Phys. Rev. A 90, 033628 (2014).
	
	\bibitem{Buecker_13} R. Bücker, T. Berrada, S. van Frank, J.-F. Schaff, T. Schumm, J. Schmiedmayer, G. Jäger, J. Grond, and U. Hohenester, J. Phys. B 46, 104012 (2013).
	
	\bibitem{Hohenester_07} U. Hohenester, P. K. Rekdal, A. Borzì, and J. Schmiedmayer, Phys. Rev. A 75, 023602 (2007).
	
	\bibitem{Pagano_24} A. Pagano, D. Jaschke, W. Weiss, and S. Montangero, Phys. Rev. Res. 6, 033282 (2024).
	
	\bibitem{Muller_22} M. M. Müller, R. S. Said, F. Jelezko, T. Calarco, and S. Montangero, Rep. Prog. Phys. 85, 076001 (2022).
	
	\bibitem{Hwang_25}  S. Hwang, H. Hwang, K. Kim, A. Byun, K. Kim, S. Jeong, M. Pratama Soegianto, and J. Ahn, Opt. Quantum 3, 64 (2025).
	
	
	\bibitem{Torrontegui_12} E. Torrontegui, X. Chen, M. Modugno, A. Ruschhaupt, D. Guéry-Odelin, and J. G. Muga, Phys. Rev. A 85, 033605 (2012).
	
	\bibitem{Lions_93} P. L. Lions and T. Paul, Revista Matematica Iberoamericana  9-3,  553 (1993).
	
	\bibitem{Gerard_91} P. G\'erard, Comm. Part. Diff. Equations 16, 1761 (1991).  
	
	\bibitem{Markowich_89} P. A. Markowich, Math. Meth. Appl. Sci. 11, 459 (1989).
	
	\bibitem{Markowich_Ringhofer} P. A. Markowich, C. A. Ringhofer, Z. Angew. Math. Mech.  69, 5675 (1989).
	
	\bibitem{Panati_02} G. Panati, H. Spohn,  S. Teufel, Phys. Rev. Lett. 88, 250405 (2002).
	
	
	\bibitem{ArnoldRinghofer1996}  A. Arnold, C. Ringhofer, SIAM journal on numerical analysis 33, 1622  (1996).
	
	\bibitem{Pulvirenti_06} M. Pulvirenti, J. Math. Phys. 47, 052103 (2006).
	
	\bibitem{Gat_14} O. Gat, M. Lein and S. Teufel, Ann. Henri Poincaré 15, 1967 (2014), 
	
	\bibitem{Denittis_03} G. De Nittis and M. Lein, Ann. Henri Poincaré 18 (2017), 1789 (2017). 
	
	\bibitem{Wu_12}	H. Wu, Z. Huang, S. Jin, D. Yin, Comm. Math. Sci., 10(4),  1301 (2012).
	
	\bibitem{Sparber_03} C. Sparber, P. Markowich, N. Mauser, Asymptot. Anal. 33, 153 (2003).
	
	
	\bibitem{Markowich_10} P. Markowich, T. Paul, C. Sparber, J. Funct. Anal. 259,  1542 (2010).
	\bibitem{Markowich_12} P. Markowich, T. Paul, C. Sparber, Arch. Ration. Mech. Anal. 205, 1031 (2012).  
	
	\bibitem{Morandi_22_PLA} O. Morandi, Phys. Lett. A 443,  128223 (2022).
	
	\bibitem{Figalli_14} A. Figalli, C. Klein, P. Markowich, C. Sparber, Comm. Pure Appl. Math. 67, 581 (2014).
	
	
	\bibitem{Antonica_17}  N. Antonica, M. Ercega, M. Lazarb, J. of Funct. Anal. 272, 3410 (2017).
	
	\bibitem{Morandi_PRA_25} O. Morandi, S. Nicoletti,  V. Gavryusev and L. Fallani, Phys. Rev. A 111, 063312 (2025). 
	
	\bibitem{Morandi_SIAM_25} O. Morandi, N. Rotundo, A. Borzì, and L. Barletti, SIAM J. Appl. Math. 84, 387 (2024).
	
	\bibitem{Manchon_15} A. Manchon, H. C. Koo, J. Nitta, S. M. Frolov and R. A. Duine, Nature Mater 14, 871 (2015).
	
	\bibitem{Morandi_25_KRM} O. Morandi, Kinetic and Related Models, 18, 76 (2025).  
	
	\bibitem{Morandi_24_JPA} O Morandi, J. Phys. A: Math. Theor. 57 145202 (2024).
	
	
	
	\bibitem{Cavaliere_19} 	F. Cavaliere, N. Traverso Ziani, F. Dolcini, M. Sassetti, F. Rossi, B 100, 155306 (2019).
	
	\bibitem{Reed_Simon_vol_II} M. Reed, and B. Simon, \emph{Functional Analysis, Methods of Modern Mathematical Physics. Vol. 2}, Academic Press, INC., London. (1980).
	
	
	
	
	
	
	%
	
	
	
	
	
	
	
	
	
	
\end{thebibliography}

							


						\end{document}